\documentclass[times,twocolumn,final,authoryear]{elsarticle}

\usepackage{framed,multirow}

\usepackage{amssymb}
\usepackage{latexsym}

\usepackage{amsmath}
\usepackage{bm}
\usepackage{makecell} 

\usepackage{url}
\usepackage{xcolor}
\definecolor{newcolor}{rgb}{.8,.349,.1}

\usepackage[citebordercolor=white]{hyperref}

\journal{Advances in Space Research}

\begin{document}


\begin{frontmatter}

\title{Manoeuvre detection in Low Earth Orbit with Radar Data}%

\author[1]{Jose M. Montilla}
\ead{jmontillag@us.es}
\author[1]{Julio C. Sanchez}
\ead{jsanchezm@us.es}
\author[1]{Rafael Vazquez}
\ead{rvazquez1@us.es}
\author[1]{Jorge Galan-Vioque}
\ead{jgv@us.es}
\author[2]{Javier Rey Benayas}
\ead{jreyb@indra.es}

\author[3]{Jan Siminski}
\ead{Jan.Siminski@esa.int}


\address[1]{Universidad de Sevilla, Escuela T\'ecnica Superior de Ingenier\'ia, Camino de los Descubrimientos s.n 41092 Sevilla, Spain}
\address[2]{Indra Sistemas, Crta Loeches 9, 28850 Torrej\'on de Ardoz, Madrid, Spain}
\address[3]{Space Debris Office, ESA/ESOC, Darmstadt, Germany}


\begin{abstract}
This work outlines and assesses several methods for the detection of manoeuvres in Low Earth Orbit  (LEO) from surveillance radar data. To be able to detect manoeuvres, the main starting assumption is that the object under analysis has an orbit known with a sufficient degree of precision. Based on the precise (a posteriori) orbit and radar data, several manoeuvre detection methods are presented; one is based on unscented Kalman filtering, whereas two others algorithms are based on reachability analysis of the state, which correlates its prediction set with the next track from the radar. The filtering algorithm can be extended for several radar tracks, whereas the reachability-based methods are more precise in detecting manoeuvres. Then, to inherit the best properties of both classes of algorithms, a manoeuvre detection filter that combines both concepts is finally presented. Manoeuvre detection results are presented first for simulated scenarios---for validation and calibration purposes---and later for real data. Radar information comes from the Spanish Space Surveillance Radar (S3TSR), with real manoeuvre information and high-quality ephemerides. The results show promise, taking into account that a single surveillance radar is the only source of data, obtaining manoeuvre detection rates of more than 50\% and false positive rates of less than 10\%.

\end{abstract}


\end{frontmatter}


\section{Introduction}
\label{sec1}

In the field of Space Surveillance and Tracking (SST), accurate orbital determination and manoeuvre detection is of upmost importance to infer objects orbital information and their future behaviour, as well as to be able to carry out tasks such as prediction of potential conjunctions with operating satellites, taking avoidance orbital corrections, predicting re-entries, identifying fragmentations or updating orbital elements of known satellites, among others.

Satellites performing unknown manoeuvres pose a challenge when trying to associate the new collected observations (obtained by laser, radar, or by any other means from the SST infrastructure) with the previously known reference orbits (which are stored in SST catalogues). Indeed, one of the main motivations of manoeuvre detection is that it can significantly reduce the number of uncorrelated targets detected by the SST sensors infrastructure. Most of these uncorrelated objects are known satellites, which have performed unpublished manoeuvres, in such a way that their new orbits do not match with the predictions.

This work develops several methods for the detection of manoeuvres in Low Earth Orbit (LEO) from radar data, providing first some preliminary numerical initial results obtained from simulated orbits and radar data. Since the final aim is to integrate these algorithms in the S3T (Spanish Space Surveillance Tracking) Cataloguing System in order to provide routine automatic manoeuvre detection capabilities to the system in the future, a validation of all the algorithms is carried out with real tracks from S3TSR~\citep{gomez2019initial}, the Spanish surveillance radar developed, installed and validated by Indra with the funding of the Spanish Government under the technical and contractual management of ESA on behalf of  Centro de Desarrollo Tecnol\'ogico e Industrial (CDTI). Manoeuvre information and ephemerides are obtained from ESA/ESOC and DLR/GSOC to assess the results, for several scenarios. The methods were implemented by making use of the space-dynamics library~\cite{orekit}. 

The contribution of this paper, compared with other existing works in the literature \citep{pastor2020satellite}, is the development of algorithms and metrics that can work with a single surveillance radar, providing a rather small number of measurements (tracks are usually 5-20 plots), large gaps in between tracks (1 to 3 days), in Low Earth Orbit, and for moderately small manoeuvres. In addition, new manoeuvre probability metrics are derived and tested, not only in simulation, but also with real data. Part of these results were already presented in conference form in~\cite{vazdebris2021,vazstardust2021}. 

The structure of this paper is as follows. After this brief introduction, a literature review is performed for the two main families of methods that can be used to detect manoeuvres, namely: Kalman filters (based on orbit determination approaches) in Section \ref{sec:two}, and reachability analysis-based methods (which compare reachable predicted sets with obtained measurements) in Section \ref{sec:three}. The particular implementations selected for this work are presented in Sections \ref{sec:four} and \ref{sec:five}, respectively for each family, together with some preliminary proof-of-concept results. Next validation results are presented for simulated (Section \ref{sec:six}) and real scenarios (Section \ref{sec_seven}). Then, in Section \ref{sec:MDF} a manoeuvre detection filter based on a combination of filtering and reachability is developed and its results presented. The paper is concluded in Section \ref{sec:eight} with some final remarks and future work.

\section{Manoeuvre detection filters}
\label{sec:two}
Manoeuvre detection filters (MDFs) employ orbit determination in the process of detecting if some manoeuvre has been performed; they are quite useful, since they are able to correlate new (post-manoeuvre) orbits with previous (pre-manoeuvre) known orbits, thus paving the way to perform orbit determination using quite fewer measures (as compared to a conventional orbit determination problem). This fact is shown, for instance, in \cite{goff2015orbit} where the author compares the accuracy and cost of orbit determination using a manoeuvre detection filter on known flying objects, versus a conventional (cold-started) orbit determination procedure. At its simplest, a manoeuvre detection algorithm relies on Statistical Orbit Determination (SOD) methods (see \cite{schutz2004statistical} for a general overview).  Very much has been written about the estimation and tracking of spacecraft using radar (or laser) measures, and classical methods like the Batch Least Squares (BLS) method, Extended Kalman Filter (EKF), and the Unscented Kalman Filter (UKF)---and even non-gaussian techniques like the Particle Filter---are well known in the literature.
 
The problem with traditional methods arises when the target performs unknown manoeuvres in-between the measurements windows. Then, the propagated orbit in which the estimation methods are based may become too inaccurate (since they do not take into account the manoeuvres) as they have become ``overconfident'' due to their covariance becoming too small (a filter exhibiting such behaviour is known as a ``smug'' filter). Thus, a manoeuvring scenario may lead to severe outliers and convergence problems in conventional filtering techniques; there exist methods to handle these issues and enhance the robustness of these classical algorithms, to avoid divergence problems. For instance, some possible techniques are covariance inflation or fading memory, among others. These mechanisms would be activated if a manoeuvre is detected; using the filter residuals, a decision logic can be put in place to estimate when a manoeuvre has been performed. 

Although not considered for this work there exist solutions with extra layers of complexity for problems with highly unstructured uncertainty. These are the Multiple Model filters, suited for tracking problems and based on a family of elemental filters which can be designed to model different aspects of the system behaviour, together with a probabilistic mixing logic to select the best estimation combining all the outputs of the elemental filters \citep{li2005survey}.

\subsection{State of the art on manoeuvre detection filters}

Next, a representative sample of MDFs from the literature are analysed.

\cite{woodburn2003estimation} presents a fixed interval smoother for manoeuvre reconstruction; this algorithm gives a simple approach that consists of a sequential filter used to move forward across the manoeuvre and a fixed interval smoother to move backwards across the manoeuvre. The sequential filter serves to process all the tracking data prior to the manoeuvre to provide an optimal pre-manoeuvre state estimate and covariance. Radar data is processed after the time of the manoeuvre until the uncertainty in the state estimate returns to a normal non-manoeuvre condition. The difference between the post-manoeuvre and pre-manoeuvre smoothed states may now be extracted as the estimate of the impulse (with covariance information used to characterize manoeuvre uncertainty). As drawback of this method it is assumed, without guarantees, that the filter will quickly converge on a good state estimate after a manoeuvre has passed.

A variable structure estimator is proposed in \cite{guang2018non}, where a manoeuvre detection metric is used to design an estimator with an additional manoeuvre observer module, the so-called ``variable structure estimator''. In this scheme an EKF is used together with a manoeuvre observer (which is in turn triggered when the manoeuvre detection metric reaches a certain threshold). Then, the manoeuvre observer estimates the manoeuvre acceleration, and sends that information to the EKF, which takes into account the estimated acceleration to improve the orbit propagation in its algorithm. The ``manoeuvre observer'' is based on a simple first-order observer, which produces an estimation of the acceleration to be fed back to the EKF. The method is simple, but assumes that  radar measurements are always available (with a frequency of 5 Hz) so that manoeuvres are always observed; that would require having data from a very large radar network.

A joint kinematic/dynamic filter is proposed in \cite{ye2021maneuver} (with simulations of three radars and a satellite in a Medium Earth Orbit) and consists of two filters running in parallel. The first and main filter is a traditional orbit determination Kalman filter whereas the second filter is a kinematic filter and utilizes some representative random processes (with design parameters) to describe the orbital motion. While the detailed motion is not captured at all, the changes caused by orbital manoeuvres can be captured by those flexible random processes, although it requires a complex tuning. 

A Filter-through and manoeuvre reconstruction approach is presented in \cite{goff2015parameter} (see also \cite{goff2015tracking}). The results  therein show that a filter-through Interacting Multiple Model orbit determination filter (EKF or UKF) can converge on a post-manoeuvre orbit with similar performance to Initial Orbit Determination (IOD) approaches, based on multiple filters running with different levels of covariance inflation. Once the post-manoeuvre orbit is known with a certain degree of accuracy, to reconstruct a single manoeuvre, one determine the time when the orbits intersected or became very close. 

\section{Reachability analysis}
\label{sec:three}
In this section the idea of Reachable Sets (RS) and their analysis (Reachability Analysis, RA) is introduced, as well as their relationship with control and estimation of systems. The technique has been used, for instance, in the context of rendezvous of spacecraft \citep{sanchez2019event}, and has many applications in the area of safety for trajectories of vehicles \citep{xu2019collision}. This background material constitutes the foundational framework for Section \ref{sec:five}.

Citing \cite{holzinger2009reachability}, ``the concept of reachability is central to Space Situational Awareness (SSA),'' which underscores the interest of this concept for the present work. Reachability Analysis deals with the study and applications of Reachable Sets, which are defined as follows: given a system that evolves from an initial condition (or set of initial conditions), and possibly has some control inputs, the reachable set is the set of states at which the system can arrive (i.e., the states that can be reached) at a given time.

To more formally define a RS, let us consider the system governed by the differential equation  $\dot{x}=f(t,x,u) $ (where $ x $ is the state, $t$ the time and $u$ a possible control input), which for a given initial condition and control spawns a trajectory $ x(t) $. This solution, if the dependencies are explicit, defines the state trajectories flow $ x(t) = \varphi^{t,t_0} (x_0,u) $. Consider an initial set (instead of a point) of initial conditions at $ t_0 $, and denote it by $ \Omega_0 $. Consider the set of all possible actuations $ U $. Then, the RS from $ \Omega_0 $ at time $ t $, denoted as $ \Omega(t) $, is defined (assuming there are no collisions or singularities for the flow) as

\begin{equation}
	\Omega (t) =  \{ x \in \mathbb{R}^n : x = \varphi^{t,t_0} (x_0,u),x_0 \in \Omega_0,u \in U \}.
\label{eq:first}
\end{equation}

Even in the linear case, the dependence on the control and initial conditions can make the computation of these sets quite difficult. If the dynamics are non-linear (as in orbital mechanics), a state transition matrix is not available, and therefore the computation of reachable sets becomes highly intensive \citep{kurzhanski2000ellipsoidal}. Our approach interprets the differential equation in a stochastic sense \citep{jain2019stochastic}, with initial conditions given as a certain initial probability distribution, so that one can consider the starting set $ \Omega_0 $ a confidence region of that probability distribution, and the RS its evolution, $ \Omega(t) $. In principle, with six states (three pertaining to position and three to velocity) that may have some degree of uncertainty, one would require to propagate the boundary of a six-dimensional closed manifold, as well as the probability distribution function inside of it. 

Thus, in this project a particle-based approach is applied instead (very much in the spirit of the Montecarlo method), in which one samples the initial confidence region, to then propagate those sample points. Since a large number of particles (trajectories) need to be propagated, the use of differential algebra techniques such as Taylor expansion over an initial condition can be employed in order to obtain $ \Omega(t) $ in a reasonable amount of time (see \cite{armellin2010asteroid} and \cite{perez2013jet}). Notice that this in fact represents a higher-order approach than the classical propagation of covariances (linear approach), that rely on Jacobians and the assumption that Gaussian distributions keep being Gaussian, which does not hold true here since the non-linearity of the propagation ``distorts'' the distribution, making it lose its Gaussian shape \citep{holzinger2009reachability}. 

In addition, one of the most interesting applications to SST of reachability analysis is the problem of object correlation. Looking at the literature, this problem has indeed received considerable attention in the last years. There are a number of metrics that can be used such as the Mahalanobis distance \citep{hall2019probabilistic} and techniques that can help when several measurements are present, such as the use of attributables \citep{vananti2017tracklet,reihs2021application}, but these do not explicitly take into account the possibility of manoeuvring objects, which is critical since small orbital corrections can produce outsized state discrepancies at the long term. This problem is tackled in \cite{singh2012space,holzinger2012object} computing (by means of optimal control) the minimum possible manoeuvre that connects the previous orbit with the new measurements. In \cite{siminski2017assessment} this optimal control approach is compared with the use of historical data, which is found more accurate when available (at least for the GEO example considered in that paper) and if the manoeuvres are predictable. These ideas are used in this work to develop manoeuvre detectors.

\section{An orbit determination filter with basic manoeuvre detection capabilities}
\label{sec:four}
To decide which filter to develop for this work, it is important to take into account that the the Spanish survey radar S3TSR \citep{gomez2019initial} is the only source of measurement data considered for this project. Being a single radar, this implies that objects will have long windows without observation in-between, from about half a day up to 3 days, and then a radar track, typically with 5-20 individual plots, will become available. Therefore, designs relying on a large number of measurements and/or frequent measurements are not implementable. The scheme of \cite{guang2018non} is adapted, with manoeuvre detection based on residue analysis. As for the choice of the filter type itself, the UKF seems to be the superior choice. The rationale of this choice is as follows. Since an EKF relies on linearisation to obtain the evolution of the state error covariance, scenarios with long propagation times such as the ones considered in this work may degrade its accuracy, depending on the starting covariance. To overcome this drawback, an UKF is considered instead, since it provides a higher-order approximation for covariance evolution which can withstand longer propagations. The UKF is based on the ``unscented transformation'' first proposed by \cite{julier1997new} and later improved by \cite{wan2000unscented} to compute the first two moments of the probability density distribution of a random variable given by some transformation $ y = h(x) $, assuming that the mean and the covariance of the variable $ x $ are known. The idea behind the unscented transformation is to use a set of points $ x^i $ (sigma-points) in such a way that their mapping $ y^i = h(x^i) $ can be used to accurately approximate the exact mean and covariance of $ y $ (by using a predefined set of weights).

\subsection{UKF algorithm}\label{sect-ukfalg}

Considering a system with $n$ states, given by the following process and observation equations
\begin{equation}
	\label{eq:second}
	\dot{\bm{x}}=f\left(\bm{x},t\right),
\end{equation}
\begin{equation}
	\label{eq:third}
	\bm{y}=h\left(\bm{x},t\right), 
\end{equation}
and a set of weights to estimate the mean and the covariance (denoted by $ w_m^{\left(j\right)} $ and $ w_c^{\left(j\right)} $ respectively, for $ j = 1,\ldots,2n+1 $), together with a tuning parameter $ \kappa $ (see \cite{goff2015orbit} for a description of the weights and the tuning parameter values), the UKF algorithm is (obtained from \cite{goff2015orbit}):

\begin{enumerate}
	\item 	Start from the previous estimate of the state and the covariance of its error ($ {\hat{\bm{x}}}_0 $ and $ \hat{\bm{\mathrm{P}}}_0 $ at first).
	\item 	Read the next observation and its covariance: $t_i,\bm{y}_i,\bm{\mathrm{R}}_i$.
	\item 	Perform the decomposition $\bm{\mathrm{P}}_{i-1}=\bm{\mathrm{A}}^T\bm{\mathrm{A}}$, and denote $\bm{a}^{\left(j\right)}$ as the \textit{j}-th column of $\bm{\mathrm{A}}$.
	\item 	Calculate the sigma points:
	\begin{equation*}
		\tilde{\bm{x}}_{i-1}^{\left(j\right)} = \hat{\bm{x}}_{i-1} + \breve{\bm{x}}^{\left(j\right)},\ \ \mathrm{for} \ \ j \ = \ 0,\dots,2n, \ \ \mathrm{and} \ \ \breve{\bm{x}}^{\left(0\right)}=0,
	\end{equation*}
	\begin{equation*}
		\breve{\bm{x}}^{\left(j\right)} = \bm{a}^{\left(j\right)} \sqrt{n +
		\kappa}, \ \ \mathrm{for} \ \ j \ = \ 1,\dots,n,
	\end{equation*}
	\begin{equation*}
		 \breve{\bm{x}}^{\left(n+j\right)} = -\bm{a}^{\left(j\right)} \sqrt{n +
			\kappa}, \ \ \mathrm{for} \ \ j \ = \ 1,\dots,n.
	\end{equation*}
	\item 	Propagate all the sigma points using numerical integration: Initial conditions $\tilde{\bm{x}}_{i-1}^{\left(j\right)}$, differential equation $\dot{\bm{x}}=f\left(\bm{x},t\right)$, integration results $\tilde{\bm{x}}_{i}^{\left(j\right)}$.
	\item Calculate the propagated a priori state and covariance (adding the process noise as in Section \ref{sec:PNoise}):
	\begin{equation*}
		\bar{\bm{x}}_i = \Sigma^{2n}_{j=0} w_m^{\left(j\right)} \tilde{\bm{x}}_{i}^{\left(j\right)}, \ \bm{\bar{\mathrm{P}}}_i = \Sigma^{2n}_{j=0} w_c^{\left(j\right)} \left(\tilde{\bm{x}}_{i}^{\left(j\right)} - \bar{\bm{x}}_{i} \right)\left(\tilde{\bm{x}}_{i}^{\left(j\right)} - \bar{\bm{x}}_{i} \right)^T + \bm{\mathrm{Q}}_i.
	\end{equation*}
	\item 	Transform the sigma-points and calculate the predicted observation:
	\begin{equation*}
		\tilde{\bm{y}}_{i}^{\left(j\right)} = h(\tilde{\bm{x}}_{i}^{\left(j\right)},t_i), \ \bar{\bm{y}}_{i} = \Sigma^{2n}_{j=0} w_m^{\left(j\right)} \tilde{\bm{y}}_{i}^{\left(j\right)}.
	\end{equation*}
	\item Calculate the predicted observation covariance and results:
	\begin{equation*}
		\bm{\mathrm{S}}_i = \Sigma^{2n}_{j=0} w_c^{\left(j\right)} \left(\tilde{\bm{y}}_{i}^{\left(j\right)} - \bar{\bm{y}}_{i} \right)\left(\tilde{\bm{y}}_{i}^{\left(j\right)} - \bar{\bm{y}}_{i} \right)^T + \bm{\mathrm{R}}_i, \bm{\nu}_i=\bm{y}_{i} -\bar{\bm{y}}_{i}.
	\end{equation*}
	\item 	Compute the Kalman gain and update the estimate of the state and its covariance:
	\begin{equation*}
		\bm{\mathrm{V}}_i = \Sigma^{2n}_{j=0} w_c^{\left(j\right)} \left(\tilde{\bm{x}}_{i}^{\left(j\right)} - \bar{\bm{x}}_{i} 	\right)\left(\tilde{\bm{y}}_{i}^{\left(j\right)} - \bar{\bm{y}}_{i} \right)^T, \ \bm{\mathrm{K}}_i = \bm{\mathrm{V}}_i (\bm{\mathrm{S}}_i)^{-1},
	\end{equation*}
	\begin{equation*}
		{\hat{\bm{x}}}_i = {\bar{\bm{x}}}_i + \bm{\mathrm{K}}_i \bm{\nu}_i, \ \hat{\bm{\mathrm{P}}}_i = \bar{\bm{\mathrm{P}}}_i - \bm{\mathrm{K}}_i \bm{\mathrm{S}}_i \bm{\mathrm{K}}_i^T.
	\end{equation*}
	\item 	Return to step 1 and process the next observation.
\end{enumerate}

\subsection{Smoothing}

In the simulated scenario considered in this work, measurements from the radar come in tracks of 5-20 plots, with long intervals in-between them (hours).  While the BLS approach is simultaneous in nature, the KF approaches (EKF/UKF)  process the measurements sequentially, in the order they were obtained; thus, the output of the filter can be improved via a backwards smoother. This additional algorithm propagates the filter backwards in time, starting from the last plot in a track, up to the first one (or even to previous tracks), modifying estimates accordingly \citep{goff2015orbit}. It is well-known that smoothers provide considerable improvement for orbit determination.

\subsection{UKF measurements}

As a first step, the UKF must be tuned to work correctly in the absence of manoeuvres. The radar measurements are range, range rate, azimuth and elevation. The radar accuracy on range is in the order of metres, whereas one gets errors under a metre per second for range rate, but the angular error translates into a distance error of kilometres for a LEO satellite, as typical for a surveillance radar of this size. Thus, only range and range rate measurements are considered as the filter's inputs.

\subsection{UKF process noise estimation}
\label{sec:PNoise}
The UKF algorithm requires the process noise covariance as an input. This quantifies mismatches with respect to the real process. Consequently, it is a key factor in the filter as it will balance the credibility of the process with respect to the measurements. In any case, the process covariance is unknown, as its exact knowledge would imply perfect modelling, and has to be tuned. Initial covariance needs also to be estimated to be as realistic as possible \citep{poore2016covariance}. In \cite{carpenter2018navigation}, the state noise compensation technique described next is recommended as a good practice for navigation filters and has been adopted.

Denote by LVLH a Local-Vertical, Local-Horizontal frame. Assuming LVLH velocity error as Gaussian white noise with covariance
\begin{equation} \label{eq:five}
	\bm{\mathcal{S}}_{LVLH} = \left[ \begin{array}{ccc} q_x & 0 & 0 \\ 0 & q_y & 0 \\ 0 & 0& q_z \end{array} \right] ,
\end{equation}

then, the transformation to inertial coordinates can be made using the rotation matrix $ \bm{\mathrm{R}}_{LVLH} $  which transforms LVLH coordinates to the inertial frame
\begin{equation} \label{eq:six}
	\bm{\mathcal{S}}_k=\bm{\mathrm{R}}_{LVLH}\bm{\mathcal{S}}_{LVLH}\bm{\mathrm{R}}_{LVLH}^T.
\end{equation}

Dividing the elapsed time between radar tracks in increments $ \Delta t $, where the inertial orientation of the LVLH frame is assumed constant, the full state inertial covariance grows during an interval $ k $ as a second-order random walk model:
\begin{equation} \label{eq:seven}
	\bm{\mathrm{M}}_k=  \left[ \begin{array}{cc} \bm{\mathcal{S}}_k \Delta t^3/3 & \bm{\mathcal{S}}_k \Delta t^2 / 2 \\ \bm{\mathcal{S}}_k \Delta t^2 / 2 & \bm{\mathcal{S}}_k \Delta t \end{array} \right],
\end{equation}

and then process noise covariance estimation at the time instant $i$ of the UKF algorithm is
\begin{eqnarray*}
	\bm{\mathrm{Q}}_i=\left\{ \begin{array}{ll} \sum_{k=1}^{N} \bm{\mathrm{M}}_k,& \textnormal{for filter calls between tracks}, \\
	0_{6\times6},& \textnormal{for filter calls within a track.} \end{array} \right.
\end{eqnarray*}

Where $\mathrm{N}$ is the number of time increments in-between filter calls. Within a track and between plots, where measurements are obtained every few seconds, the process mismatch is negligible. 

\subsection{UKF preliminary testing results}
\label{sec:four_six}
\begin{figure}[t!]
		\centering
	\includegraphics[scale=0.44]{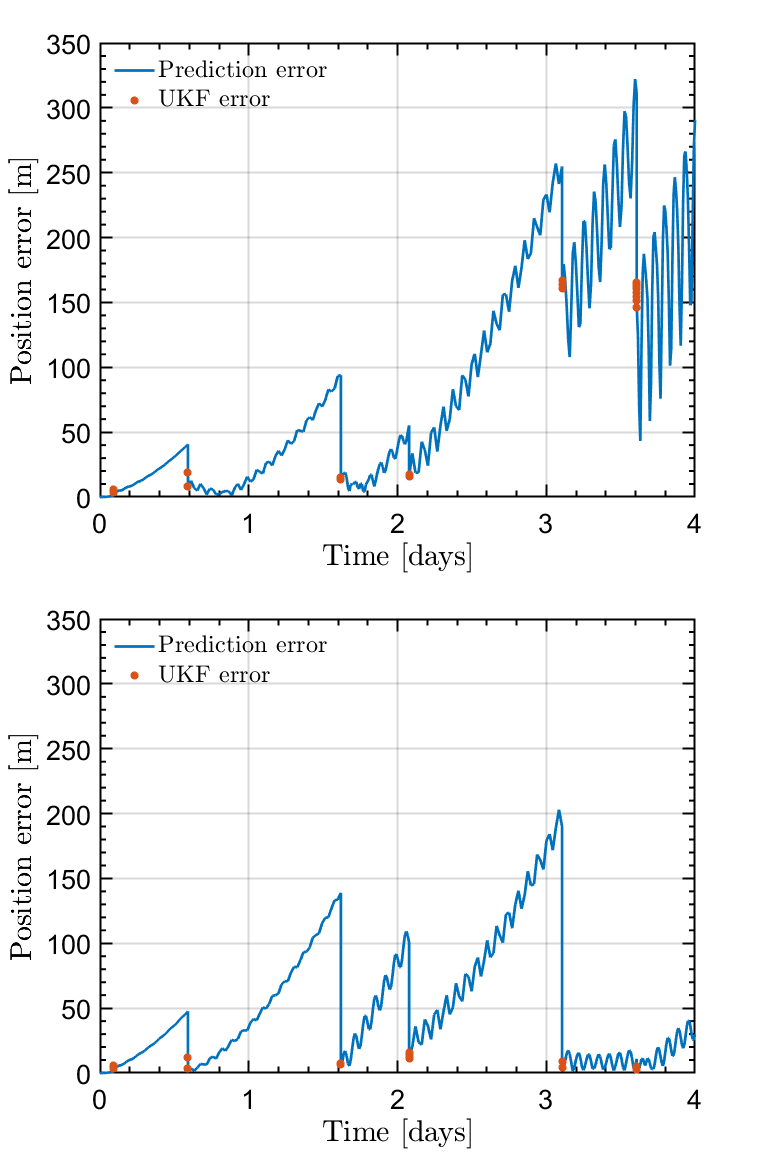}
	\caption{Position error without process noise (top) and estimating process noise (bottom). Red dots indicate the UKF error at each plot (measurement points).}
	\label{fig:first}
\end{figure}

Numerical results are shown in order to justify the chosen implementation. The considered scenario is the LEO satellite Sentinel-1A (with initial orbital elements taken from public TLEs and assumed precise, propagated with second-order gravity harmonics and drag using OREKIT) between 16:00:00 08/07/2015 to 16:00:00 12/07/2015. The radar measurements are simulated with some noise. The following results assume a model mismatch in drag, with $ C_D=2.2 $, $ S=10\ \mathrm{m}^2 $ the ``real'' drag coefficient and exposed surface, and $ C_D=2 $, $ S=9.5\ \mathrm{m}^2 $ the assumed ones. The LVLH acceleration errors for (\ref{eq:five}) are
\begin{eqnarray*}
	q_x={10}^{-9}\ \mathrm{m}^2s^{-3},\ q_y=q_z=5 \cdot 10^{-10}\ \mathrm{m}^2s^{-3}
\end{eqnarray*}
where more process noise has been assumed in the tangential direction (due to the drag modelling uncertainty being dominant in LEO). The discretization time period is taken as $ \Delta t=10\ \mathrm{min} $. A comparison with a simulation assuming a null process covariance noise is shown in Figure \ref{fig:first} (red dots indicate the mismatch between measurements and the predicted state after the filter update). It can be seen how the inclusion of some process noise greatly improves the filter's convergence.

In the results, the initial covariance is assumed small and with a realistic shape (obtained from running the filter for a few days); using a diagonal shape resulted in a much poorer performance of the filter.

\subsection{Manoeuvre detection}

The filter is adapted to estimate the presence of manoeuvres. Thus, in the 8th step of the UKF's algorithm, a manoeuvre prediction metric can be included, which reads:
\begin{equation} \label{eq:four}
	\mathrm{\Psi}_i = \sqrt{\bm{\nu}_i^T \bm{\mathrm{S}}_i^{-1} \bm{\nu}_i},
\end{equation}

with $ \bm{\nu}_i $ being the residuals and $ \bm{\mathrm{S}}_i $ the observation covariance, for $i\in\{1,\hdots,n\}$, $n$ being the number of plots of a given track. This term can be used to estimate model mismatches (due to manoeuvres), and then trigger other manoeuvre detection algorithms, every time a radar track is processed.

For manoeuvre detection, three metrics (\ref{eq:four}) were proposed; either $\Psi_1$ (with the logic that the first plot may show the largest impact from a previous manoeuvre), $\max_{i} \{\Psi_i\}$ or
\begin{equation}
	\Psi = \frac{1}{n} \sqrt{\sum_{i=1}^n \Psi_i },
\end{equation}
As the detection method is based on finding significant discrepancies between the predicted orbit and the actual one, the residuals are computed prior to smoothing. 
\begin{figure}[t!]
	\centering
	\includegraphics[scale=0.44]{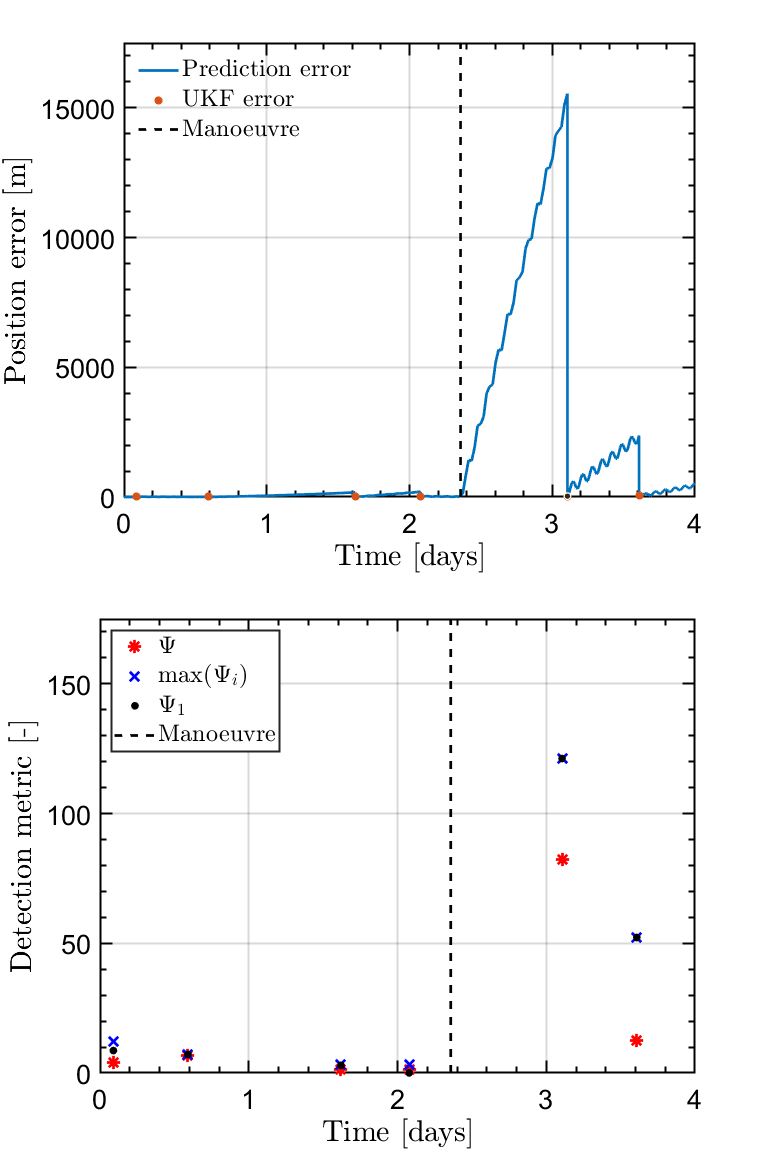}
	\caption{UKF position error with manoeuvre (top) and detection metrics comparison (bottom). The vertical dashed line represents the manoeuvre}
	\label{fig:two}
\end{figure}
A proof of concept is shown for the same scenario with a manoeuvre. The manoeuvre start is 00:34:58 11/07/2015 and ends at 00:35:24 11/07/2015 with a constant acceleration of $ u=\left[ 0.31 \cdot 10^{-2},-0.35 \cdot 10^{-3},0.37 \cdot 10^{-5} \right] m/s^2 $ in the LVLH frame. In Figure \ref{fig:two}, the UKF demonstrates its capability to recover the orbit after the manoeuvre is applied, even without any particular mechanism for covariance inflation. Moreover, a comparison of the possible detection metrics is also shown, with $\Psi$ seemingly being the metric with a better trade-off between detection and false positives.

\section{Reachability-based manoeuvre detection algorithms}
\label{sec:five}
In this section, RA as outlined in Section \ref{sec:three} is applied to the specific problem of manoeuvre detection. Thus, the starting inputs are the precise orbit of the objects (pre-manoeuvre) and a radar track (possible, post-manoeuvre), and the output is the probability of a manoeuvre having happened. 

As a first step, the theory of attributables \citep{vananti2017tracklet,reihs2021application} is introduced; it allows to ``compress'' several plots into a single, higher-quality measurement, fitting a full track into a single polynomial expression whose order needs to be determined.

Next, two algorithms are explained; the first is based on comparing the range and range rate attributables ($\rho,\dot{\rho}$) obtained from measurements with the one obtained from the initial uncertain orbit, by means of confidence regions and the Mahalanobis distance, which is a measure of the distance between a point and a distribution \citep{hall2019probabilistic}.

The second algorithm is based on the use of optimal control theory. Following the ideas of \cite{singh2012space} and \cite{holzinger2012object} one can compute by means of stochastic optimal control a distribution of the $ \Delta V $ that connects the uncertain orbit with the measurement. This distribution can then be used to obtain the likelihood of a manoeuvre having been performed.

\subsection{Attributables}
\label{sec:cinco_uno}
We use attributables to condense the information of all plots in each track \citep{reihs2021application}. A radar provides range $\rho$, range rate $\dot \rho$, elevation $\mathrm{El}$ and azimuth $\mathrm{Az}$, that, coupled with the chosen reference epoch $t$, form the attributable
\begin{equation}
	\mathcal{A}=\left\{t,\rho,\mathrm{El},\mathrm{Az},\dot{\rho}\ \right\}.
\end{equation}

Fitting the information of the observables independently is one option, but it is possible to reduce the uncertainty of the resulting virtual measurement if one incorporates the definition of range-rate into the modelling, so that it shares the parameters with the range, as follows:
\begin{eqnarray}
	\label{eq:diez}
	\rho\left(t\right)=\rho_0+\rho_1t+\rho_2\frac{t^2}{2!}+\cdots+\rho_n\frac{t^n}{n!}\ ,
\end{eqnarray}
\begin{eqnarray}
	\label{eq:once}
	\mathrm{El}\left(t\right)={\mathrm{El}}_0+{\mathrm{El}}_1t+{\mathrm{El}}_2\frac{t^2}{2!}+\cdots+{\mathrm{El}}_n\frac{t^n}{n!}\ ,
\end{eqnarray}
\begin{eqnarray}
	\label{eq:doce}
	\mathrm{Az}\left(t\right)={\mathrm{Az}}_0+{\mathrm{Az}}_1t+{\mathrm{Az}}_2\frac{t^2}{2!}+\cdots+{\mathrm{Az}}_n\frac{t^n}{n!}\ ,\ 
\end{eqnarray}
\begin{eqnarray}
	\label{eq:trece}
	\dot{\rho}\left(t\right)=\frac{d\rho\left(t\right)}{dt}=\rho_1+\rho_2\frac{2t}{2!}+\cdots+\rho_n\frac{nt^{n-1}}{n!}.
\end{eqnarray}

In the expressions above, the origin of time would be at the middle time of the track.

This method manages to average out noise and reduce the standard deviation of the virtual measurement. Following the nomenclature in \cite{reihs2021application}, the set of equations that allows to solve the parameters in the sense of least-squares is:
\begin{eqnarray}
	\bm{m}= \left[ \begin{array}{c}
		\bm{\rho} \\ \bm{\mathrm{El}} \\ \bm{\mathrm{Az}} \\ \bm{\dot{\rho}}  
	\end{array} \right] = A_{SYS}  \bm{p} + \bm{\upsilon} = 
	\left[ \begin{array}{ccc}
		A & 0 & 0 \\
		0 & A & 0 \\
		0 & 0 & A \\
		A_{\dot{\rho}} & 0 & 0
	\end{array} \right] \bm{p} + \bm{\upsilon},
\end{eqnarray}
where $ \bm{m} $ contains the measurements of all observables in the track, $ \bm{p} $ collects the parameters $\rho_i$, $\mathrm{Az}_i$, and $\mathrm{El}_i$ that one wants to calculate, and the matrices $ A $ and $ A_{\dot{\rho}}$ have the coefficients in the formulas (\ref{eq:diez})--(\ref{eq:doce}) and (\ref{eq:trece}), respectively, evaluated at the times of the corresponding plots. The error $ \bm\upsilon $ is assumed to follow a Gaussian distribution. Then the problem to solve is posed using weighted least-squares:
\begin{equation}
	\label{eq:quince}
	\underset{p}{\min}\ \bm{\upsilon}^T W \bm{\upsilon} = \underset{p}{\min}\ (\bm{m}-A_{SYS}\bm{p})^T W (\bm{m}-A_{SYS}\bm{p}),
\end{equation}
whose solution is well-known:
\begin{equation}
	\label{eq:dseis}
	\bm{p} = (A_{SYS}^T W A_{SYS})^{-1} A_{SYS}^T W \bm{m}.
\end{equation}

The weighting matrix $ W $ is chosen to be the inverse of the covariance matrix of the measurements, $ \Sigma_{\upsilon} $, thus the attributable errors covariance matrix is:
\begin{equation}
	\Sigma_p = \left(  A_{SYS}^T \Sigma_{\upsilon}^{-1} A_{SYS} \right) ^{-1}.
\end{equation}

This allows to estimate how good a virtual measurement $ \alpha\left(t\right) $ is expected to be to be at any point of the fit, by computing
\begin{equation}
	Var\left[\alpha\left(t\right)\right]=\ \sum_{i,j}{\sigma_{ij}\frac{t^it^j}{i!j!}},
\end{equation}
where $\sigma_{ij} $ are the coefficients of $ \Sigma_p $ corresponding to the sub-matrix of each observable; these are sufficient to provide the covariance matrix of the complete attributable ($ \Sigma_\mathcal{A} $) at the epoch without further processing.

\begin{figure}[t!]
	\hspace*{-0.55cm}   
		\centering
	\includegraphics[scale=0.76]{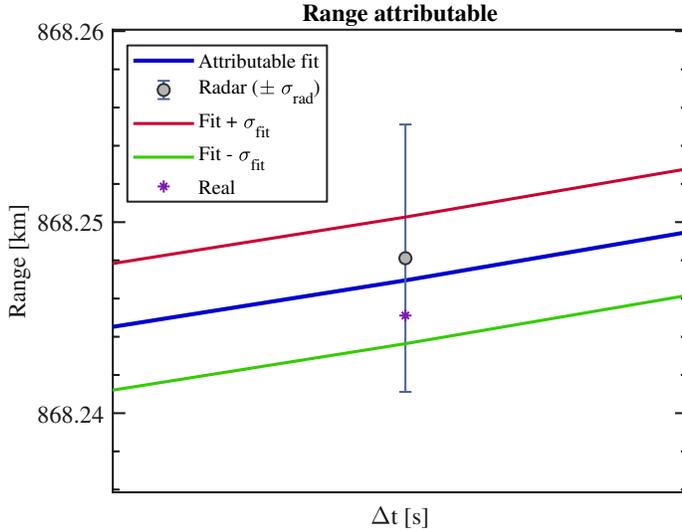}
	\caption{Range attributable and reduction of error.}
	\label{fig:three}
\end{figure}

A test track shown in Figure \ref{fig:three}, with realistic radar standard deviations, has been used as an example of range attributable. The uncertainty is mitigated with a noticeable reduction of $ 50\% $ in the standard deviation (in the case of the range rate, the reduction is of $ 60\% $), a consequence of the methodology used for the definition of the range-rate attributable.

In addition, one could even try to estimate the azimuth and elevation rates and use it for Initial Orbit Determination (IOD), but as one cannot expect to perform well as a good IOD with such a short observation arc, this method is not pursed in this work.

\subsection{Algorithm 1: Comparison of real and projected attributable through Mahalanobis distance}
\label{sec_alg1}
Applying the nomenclature of Section \ref{sec:cinco_uno}, from all the plots of a track one can obtain the virtual values of range, range-rate, azimuth and elevation at the middle of a track ($t_0$), namely $\mathcal{A}=(t_0, \rho_0,\mathrm{Az}_0,\mathrm{El}_0,\rho_1)$, as well as the associated uncertainty in the form of a covariance matrix $ \Sigma_\mathcal{A} $. 

The following algorithm is used to obtain a ``projected'' (or predicted) measurement from the initial value of the reference orbit, which is assumed to follow a certain known distribution:

\begin{enumerate}
	\item 	Sample the PDE of the initial condition obtaining $ m $ sample points. Denote these as $ x_{0j} $ for $ j=1,\ldots,m $. The set of initial conditions $ \Omega_0 $ is then approximated by these points.
	\item 	Propagate the sampled points using an OREKIT propagator up to time $ t_f $. Taylor differential algebra methods can be used to greatly speed up this computation, at the price of a lengthy initial calculation \citep{andrea2016automatic}. Thus, one obtains m trajectories $ x_j\left(t\right) $.
	\item 	Projected values at the attributable time $ t_i $ are obtained as a cloud of points $ x_j(t_i) $, with the density of points giving an approximate measure of the probability associated to the real trajectory.
	\item 	Now for each sampled orbit, one can compute the radar measurements at $t_0$, obtaining a ``cloud'' of measurements, from which one can obtain its mean ($ \widehat{\rho_0},\widehat{\mathrm{Az}_0},\widehat{\mathrm{El}_0},\widehat{\rho_1} $) as well as the associated uncertainty in the form of a covariance matrix  $ \hat{\Sigma} $. This is denoted as the projected measurement (in the sense that it is the attributable value one would expect given the  distribution of the initial condition).
	\item Finally, the attributables and projected measurements can be compared. If no manoeuvre has been performed, one would expect that both values should somewhat agree. To formulate this more precisely, define:
	\begin{eqnarray}
		\left[ \begin{array}{c}
			{\Delta\rho}_0 \\ {\Delta \mathrm{El}}_0 \\ {\Delta \mathrm{Az}}_0 \\ {\Delta\rho}_1
		\end{array} \right] = 
		\left[ \begin{array}{c}
			{\rho}_0 \\ {\mathrm{El}}_0 \\ {\mathrm{Az}}_0 \\ {\rho}_1
		\end{array} \right] - 
		\left[ \begin{array}{c}
			\widehat{\rho_0} \\ \widehat{\mathrm{El}_0} \\ \widehat{\mathrm{Az}_0} \\ \widehat{\rho_1}
		\end{array} \right], \Delta\Sigma=\Sigma_\mathcal{A}+\hat{\Sigma}.\ 
	\end{eqnarray}
	\item 	Then, if there is no manoeuvre, one would expect that, under an assumption of normality, $ \left({\Delta\rho}_0,\Delta A z_0,\Delta \mathrm{El}_0,\Delta\rho_1\right) $ should belong to a normal distribution of zero mean and covariance $ \Delta\Sigma $. This can be checked either by computing confidence regions or equivalently through the Mahalanobis distance, as briefly explained next. 
\end{enumerate}

\subsubsection{Use of confidence regions and Mahalanobis distance}

For a $n$-dimensional multivariate normal distribution with mean $m$ and covariance matrix $ \Sigma $, the $ p-\textnormal{level} $ confidence ellipsoid (this is, the ellipsoid containing with probability $ p $ samples from the distribution) is given by
\begin{eqnarray}
	\left\{ x\in \mathbb R^n:({x-m)}^T\Sigma^{-1}\left(x-m\right)\le\chi_n^2\left(p\right)\right\},
\end{eqnarray}
where $ \chi_n^2(p) $ is the inverse cumulative distribution function of the Chi-square distribution with $ n $ degrees of freedom (the dimension of the vector $ x $), evaluated at the probability value $ p $. Similarly, the Mahalanobis distance is a measure of the distance of a point $ x $ from a distribution. It is unitless, scale-invariant and takes into account the correlations of the distribution. Concretely, if the distribution has mean m and covariance matrix $ \Sigma $ the Mahalanobis distance (MD) of a point $x$ is computed as
\begin{eqnarray}
	\mathrm{\mathrm{MD}}(x)=\sqrt{\left(x-m\right)^T\Sigma^{-1}\left(x-m\right)}\ .
\end{eqnarray}

In particular if the distribution is a multivariate normal, then the MD has a chi-square distribution with $ n $ degrees of freedom; thus, it is equivalent to the use of confidence ellipsoids. This property can be used to compute probabilities of manoeuvre.
\begin{figure}[t!]
	\centering
	\includegraphics[scale=0.85]{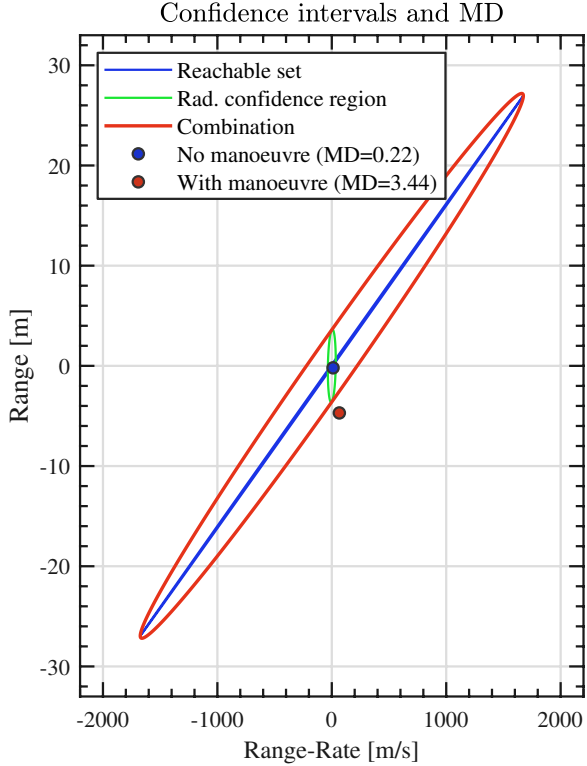}
	\caption{Confidence intervals ($ 90\% $) and Mahalanobis distance, considering only range and range-rate, without manoeuvre (blue circle) and with manoeuvre (red circle). The green ellipse represents the measurement uncertainty (radar), the blue one the orbit uncertainty, and the red one the total uncertainty.}
	\label{fig:four}
\end{figure}
Next, an example is shown where the comparison of real and projected attributables is carried out for two cases: one example with manoeuvre and one without. Figure \ref{fig:four} shows that the confidence intervals and MD are able to discriminate the manoeuvred case from the non-manoeuvred one, at least for a simple basic simulation, using range and range-rate.

Based on the MD, a probability measure has been computed based on the MD being distributed as a $ \chi^2 $ distribution function with as many degrees of freedom ($ n $) as variables, by fixing a threshold of being outside the 50\% ellipsoid. Thus, a number $\mathrm{PR}_{\mathrm{MD}}$ (which is a possible estimation of the probability of manoeuvre) is computed as follows:
\begin{equation}
	\mathrm{PR}_{	\mathrm{MD}}=max\left\{0,{2(\chi}^2\left(\mathrm{MD};n\right)-0.5\right\}. \label{eq:prmd}
\end{equation}
This way, if the MD has a probability of $ 50\% $ or less of occurring, it is assumed that there is no manoeuvre, to reduce false positives. If the MD has a probability of more than $ 50\% $ of happening, then one subtracts $ 50 $ from the probability and multiplies it by two; if one gets, e.g., a probability of a certain MD of $ 80\% $, then $ \mathrm{PR}_{	\mathrm{MD}}=60\% $.

\subsection{Algorithm 2: use of optimal control to compute a $\Delta V$ bounded measurement of distance}

As a more sophisticated alternative to the Mahalanobis distance \citep{singh2012space,holzinger2012object} one can compute by means of stochastic optimal control a distribution of the minimum $ \Delta V $ that connects the uncertain orbit around it. This distribution can then be used as a metric to obtain the likelihood of a manoeuvre having been performed. The optimal control problem is posed as follows:
\begin{eqnarray}
	\mathcal{J} &=& \min_u\ \int_{t_0}^{t_f}{u^T\left(t\right)u\left(t\right)}dt, \\
	s.t.\ \ \ \ x^\prime\left(t\right)&=&f(x\left(t\right),u\left(t\right),t) \nonumber, \\
	x\left(t_0\right)&=&x_0, \nonumber\\
	h\left(x\left(t_f\right)\right)&=&\left[\rho\ \dot{\rho}\right]^T\nonumber.
\end{eqnarray}

In the above optimal control problem, the initial point is known from the precise orbit whereas the function $h$ at the final point represents the function relating position and velocity with range and range-rate (the most precise measurements) which should take the value obtained with attributables as explained in Section \ref{sec:cinco_uno}. The function $ f $ represents the orbital dynamics, including any desired perturbation. The selected functional would represent the energy of the manoeuvre acceleration, which is less problematic than its $ \textnormal{L}_2$ norm from a numerical point of view. It is well known from the literature that the real $ \Delta V $ is bounded by the square root of this quantity, see, e.g.,~\cite{siminski2017assessment}.

The problem is solved with CasADi \citep{andersson2019casadi}, an open-source solver for MATLAB, with a multiple shooting method discretizing the orbital dynamics in N time intervals; for each of these, since manoeuvres are small, the orbital dynamics is replaced with a linearized model obtained from OREKIT (computing the State Transition Matrix, or STM), with  discrete $\Delta V$'s applied at the beginning.

\begin{figure}[t!]
	\centering
	\includegraphics[scale=0.545]{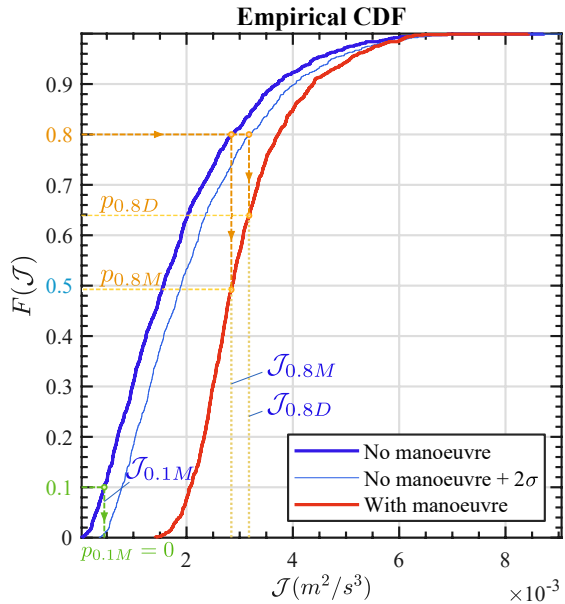}
	\caption{Empirical distribution of manoeuvre energy $ \mathcal{J} $ computed from the stochastic optimal control problem (1000 samples). The dark blue line represents the mean CDF of the non-manoeuvred case, and the light blue is the mean plus 2-sigma distribution, whereas the red plot represents a manoeuvred cases.}
	\label{fig:five}
\end{figure}

As a first step, the problem has been solved in a deterministic way. Since, once the STM is computed, the solution is fast (seconds or less), to incorporate the stochasticity of the problem (both in initial orbit and measurements), a Monte Carlo algorithm has been implemented as a simple solution, albeit rather time-consuming. Figure \ref{fig:five} shows the obtained cumulative empirical distribution of $ \mathcal{J} $ (from $ 1000 $ samples) for two cases (with and without manoeuvre).

In addition, a novel method to discriminate potential manoeuvres is now described. Qualitatively, it is clear that the distribution without manoeuvre is ``smaller'' than the one with manoeuvre. In the case without manoevre, we can derive a ``mean distribution'' as well as a distribution at a 2-sigma distance from the mean, which is helpful to avoid false positives. From these distributions some metrics have been defined, by using its $ 10\% $, $ 50\% $ and $ 80\% $ percentiles. 

The metrics are exemplified in Fig.~\ref{fig:five}. As shown, take the $ 80\% $ percentile of the estimated mean non-manoeuvred distribution, $\mathcal{J}_{0.8M}$, and compare with the probability $p_{0.8M}=\mathrm{Pr}\left(\mathcal{J}\leq\mathcal{J}_{0.8M}\right)$ of the (potentially manoeuvred) distribution to be below that energy value, as graphically shown in Fig.~\ref{fig:five}. The higher that probability, the less likely of a manoeuvre to have happened according to that particular percentile (as there would be more cases that require the same or less energy than the non-manoeuvred case to connect the initial orbit and the measurement). In the figure it can be seen how using the 2-$\sigma$ non-manoeuvred distribution (thus using $\mathcal{J}_{0.8D}$ instead, leading to $p_{0.8D}$) is more conservative. To use this idea to establish a probability of manoeuvre from a given percentile $d$ the following scaling is used:
\begin{eqnarray}
	\mathrm{P}_d=max\left\{0,\frac{\left(d-\mathrm{Pr}\left(\mathcal{J}\leq\mathcal{J}_{d}\right)\right)}{d}\right\}.
\end{eqnarray}
For instance, for $d=0.1$ and calling $p_{0.1}=\mathrm{Pr}\left(\mathcal{J}\leq\mathcal{J}_{0.1}\right)$, if $p_{0.1}$ is above $ 10\% $ the probability becomes zero, and if not, the difference is multiplied by 10, which means that $p_{0.1}=0$ would indicate total confidence of manoeuvre, for that metric (as shown in Fig.~\ref{fig:five}). When using the mean distribution we will refer to this metric as P1M, with P1D reserved for the 2-sigma metric. Other metrics we use employ the $ 50\% $ and $ 80\% $ percentiles and are denoted as P5M, P5D, P8M and P8D.

\section{Results for simulated scenarios}
\label{sec:six}
An OREKIT-based simulator, both for the manoeuvres and for the radar observations, has been developed. They provide realistic (though not accurate) testing examples. They are very useful to tune and validate the different algorithms and filters. Starting points are generated from public TLEs which are used to  define reference orbits with propagators including $\mathrm{J}_2$ and aerodynamic drag, as explained in Section \ref{sec:four_six}.

The algorithms, besides the model mismatches explained in that section, start from initial conditions within the expected limits of error of the real precise orbits (meters). Two main scenarios, respectively based on the satellites Sentinel-1A and Swarm-C, are considered. A manoeuvre either tangential (T), out-of-plane (OOP) or hybrid (with components both tangential and out of plane) is simulated, maintaining a constant acceleration of $ 10^{-3}\ m/s^2 $ and characterized by the following fields:

\begin{enumerate}
	\item Manoeuvre intensity (regulated through the duration): low (5 seg $\rightarrow 5 \cdot 10^{-3}\ m/s$) / medium (30 seg $\rightarrow 3 \cdot 10^{-2}\ m/s$) / high (120 seg $\rightarrow 1.2 \cdot 10^{-1}\ m/s$).
	\item Manoeuvre location with respect to a radar track: 2 h, 6 h or 12 h before radar.
	\item The Sentinel-1A scenario spans from 00:00:00 18/08/2020 to 00:00:00 22/08/2020. The manoeuvre starts at 18:25:00 20/08/2020.
	\item The Swarm C scenario spans from 00:00:00 14/07/2020 to 00:00:00 20/07/2020. The manoeuvre starts at 12:30:00 17/07/2020.
\end{enumerate}

Thus, combining all these factors, one gets 18 simulation scenarios per satellite to analyse the influence of these factors for the algorithms. Due to space limitations, only selected results are shown, with general conclusions drawn from the complete set.

\subsection{UKF simulated results}
\label{sec:six_one}
\begin{figure}[t!]
	\centering
	\includegraphics[scale=0.46]{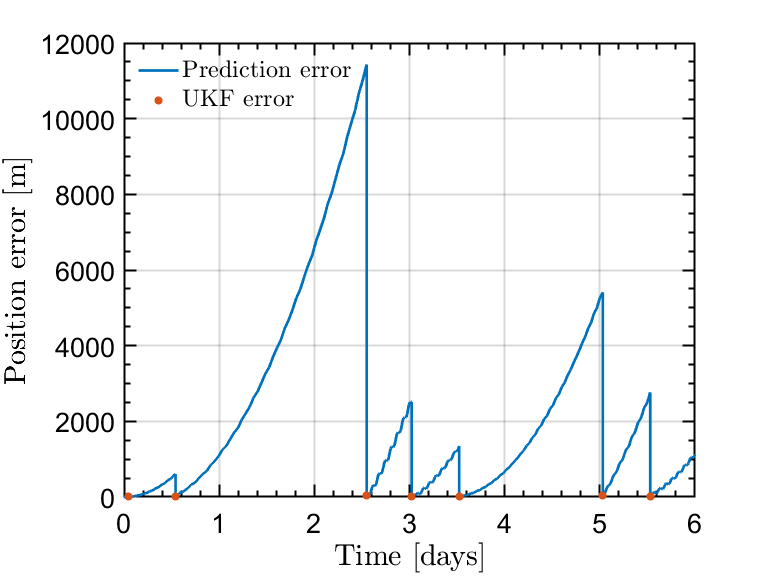}
	\caption{Position error for Sentinel-1A with respect to reference orbit without manoeuvre. Red dots indicate the mismatch between measurements and the predicted state after the filter update.}
	\label{fig:six}
\end{figure}

\begin{figure*}
	\centering
	\includegraphics[scale=0.55]{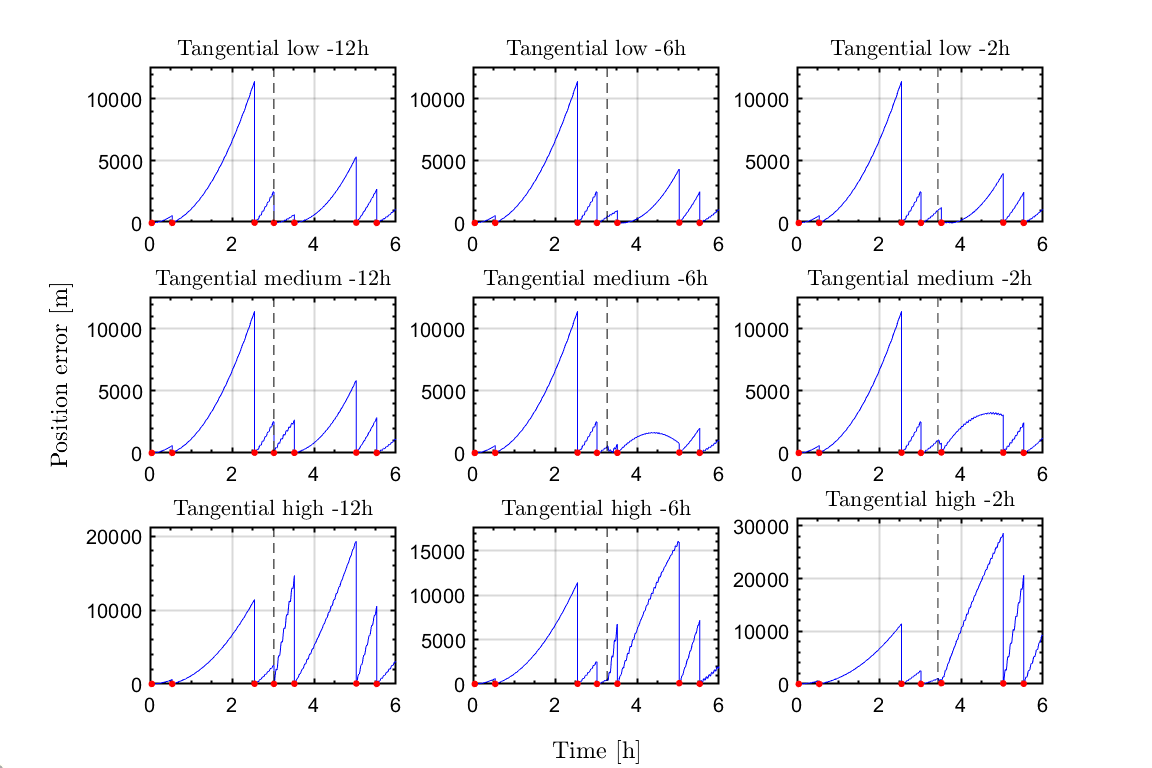}
	\caption{Position error with respect to reference orbit for Sentinel-1A manoeuvred scenarios. Red dots indicate the mismatch between measurements and the predicted state after the filter update.}
	\label{fig:eight}
\end{figure*}

The result without manoeuvre for Sentinel-1A is presented in Figure \ref{fig:six}, whereas the manoeuvred case (tangential) is shown in Figure \ref{fig:eight}. The value of $ \Psi $, which should help in detecting manoeuvres, is given for some Sentinel-1A cases in Table \ref{table:one}.

\begin{table}[h!]
	\centering
	\begin{tabular}{|l|l|l|}
		\hline
		\bm{$Case /  \Psi $}     & \bm{$ Pre-man. Max. $} & \bm{$ Post-man. Max. $} \\
		\hline
		\bm{$ No\ manoeuvre $} & 4.435        & 2.341*     \\
		\hline
		\bm{$ low -2 h $} & 4.435        & 2.386*        \\
		\hline
		\bm{$ low -6 h $} & 4.435      & 2.411*       \\
		\hline
		\bm{$ low -12 h $} & 4.435      & 2.335*       \\
		\hline
		\bm{$ medium -2 h $} & 4.435      & 3.478*       \\
		\hline
		\bm{$ medium -6 h $} & 4.435      & 10.33*       \\
		\hline
		\bm{$ medium -12 h $} & 4.435      & 6.260      \\
		\hline
		\bm{$ high -2 h $} & 4.435      & 34.39*       \\
		\hline
		\bm{$ high -6 h $} & 4.435      & 132.0*       \\
		\hline
		\bm{$ high -12 h $} & 4.435      & 75.89       \\
		\hline
	\end{tabular}
	\caption{Maximum value of filter detection metric before and after Sentinel-1A tangential manoeuvre, simulated results. The asterisk indicates that the maximum arises after the first post-manoeuvre track (i.e., at a later track)}
	\label{table:one}
\end{table}

It can be observed e.g. in Fig.~\ref{fig:six} that the filter takes some time to stabilize. This is probably due to the incorrect initial covariance. Since in real scenarios the covariance will not be perfectly known this can be expected. On the other hand, it is clear that the filter is working correctly in all cases; since the measurements are scarce it is unavoidable that the position errors grow, however, they are clearly mitigated at each measurement. It can be seen in Fig.~\ref{fig:eight} that manoeuvres induce large errors after they happen, since they are unaccounted for in the process covariance. The largest the manoeuvre the larger the error and the more it takes to recover from it. From Table~\ref{table:one}, one can observe that the value of $ \Psi $ is indicative of the presence of a manoeuvre only in medium and specially in high-intensity cases. Low-intensity manoeuvres are indistinguishable from process noise. In addition, the distance to the radar measurement does not seem to have much influence in the value of $ \Psi $

In the Swarm C case (not shown), the value of $ \Psi $ was indicative of the presence of a manoeuvre only in high-intensity cases. For low- and medium-intensity manoeuvres, they were, in principle, indistinguishable from process noise, unless the manoeuvre happened at a long enough distance from the first radar measurement. The main cause of this was, besides the long gap without measurement, having less radar measurements; in the case of Sentinel-1A, nine values were obtained at the pass after the manoeuvre, whereas in the case of Swarm-C, only five values are obtained.

\subsection{Algorithm 1  simulated results}
\label{sec:six_two}
\begin{figure*}[h!]
	\centering
	\includegraphics[scale=0.55]{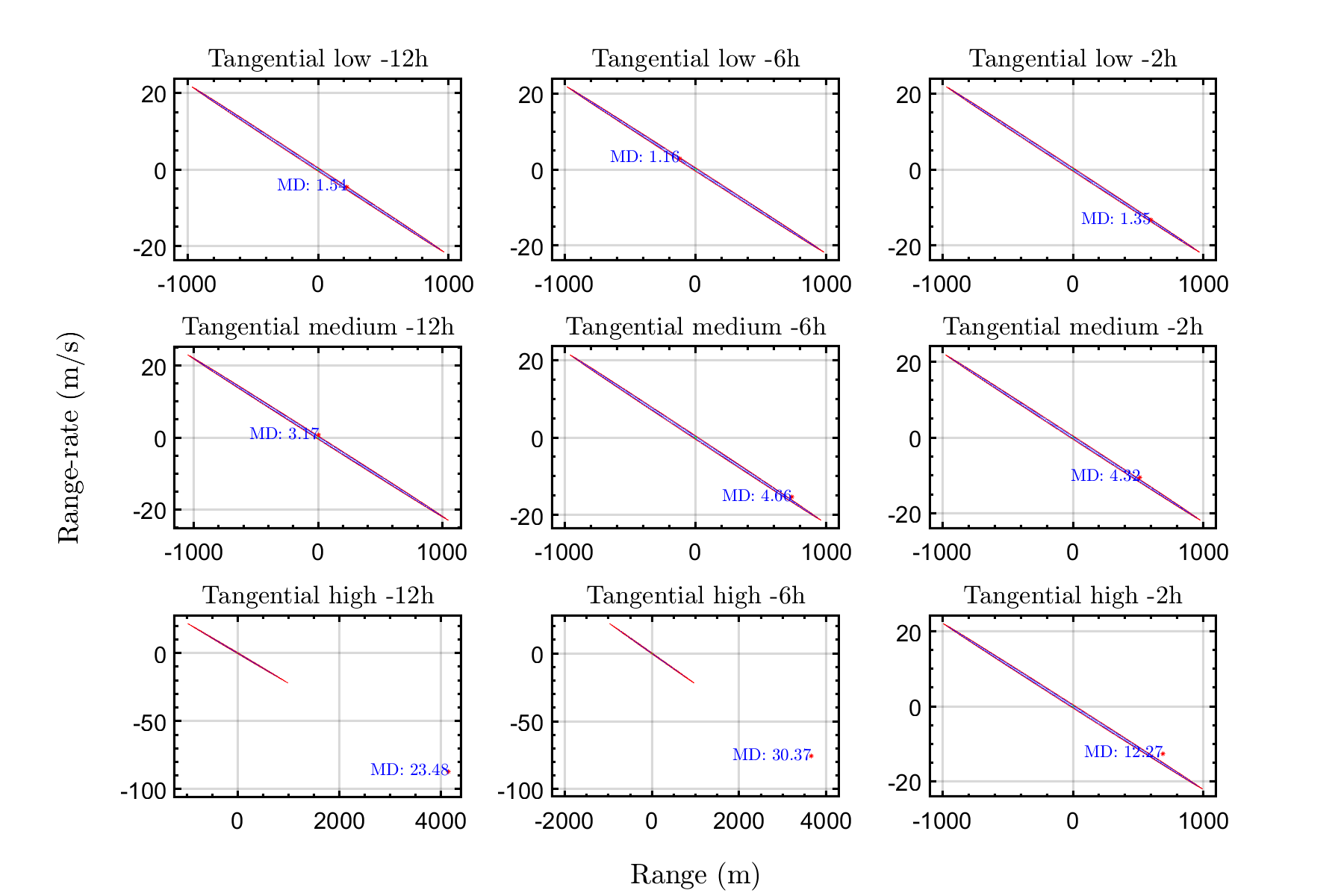}
	\caption{Reachability results (Algorithm 1) for Sentinel-1A manoeuvre, Range vs Range-rate.}
	\label{fig:nine}
\end{figure*}

In Table \ref{table:two}, the Mahalanobis distance (MD) has been computed considering only range-range rate (column 2), El-Az (column 4) and all four measurements (column 6).  As seen in the table, in general, using only range and range-rate is more sensitive in more cases; using elevation and azimuth can induce false positives. All high- and medium-intensity manoeuvres are detected, but low-intensity manoeuvres are usually not detected. In addition, the distance to the radar track does not seem to affect these results.

\begin{table*}[ht!]
	\centering
	\begin{tabular}{|l|l|l|l|l|l|l|}
		\hline
		\thead{Case / Metrics}     & \thead{MD \bm{$  (\rho,\dot{\rho}) $}} & \thead{$\mathrm{PR}_{\mathrm{MD}}$ \bm{$ (\%) $}} & \thead{MD \bm{$ (\mathrm{El},\mathrm{Az}) $}} & \thead{$\mathrm{PR}_{\mathrm{MD}}$ \bm{$ (\%) $}} & \thead{MD (All)} & \thead{$\mathrm{PR}_{\mathrm{MD}}$ \bm{$ (\%) $}} \\
		\hline
		\bm{$ WoM $} & 0.29        & 0  & 2.18 & 33 & 2.60 & 0     \\
		\hline
		\bm{$ L-12h $} & 1.54        & 8  & 1.66 & 13 & 2.22 & 0     \\
		\hline
		\bm{$ L-6h $} & 1.16        & 0  & 1.52 & 7 & 1.96 & 0     \\
		\hline
		\bm{$ L-2h $} & 1.35        & 0  & 2.10 & 30 & 2.16 & 0     \\
		\hline
		\bm{$ M-12h $} & 3.17        & 59  & 2.15 & 32 & 3.97 & 18     \\
		\hline
		\bm{$ M-6h $} & 4.66        & 81  & 0.97 & 0 & 4.81 & 39     \\
		\hline
		\bm{$ M-2h $} & 4.32        & 77  & 2.29 & 36 & 4.80 & 38     \\
		\hline
		\bm{$ H-12h $} & 23.48        & 100  & 3.61 & 67 & 23.53 & 100    \\
		\hline
		\bm{$ H-6h $} & 30.37       & 100  & 4.44 & 78 & 30.38 & 100     \\
		\hline
		\bm{$ H-2h $} & 12.27        & 100  & 1.23 & 0 & 12.44 & 97     \\
		\hline
	\end{tabular}
	\caption{Sentinel-1A reachability analysis with Algorithm 1 and probability from MD. WoM=without manoeuvre, L=low, M=medium, H=high. Tangential case}
	\label{table:two}
\end{table*}

The results can be inspected visually in Figure \ref{fig:nine}. Note that due to the propagation ``stretching'' the orbit uncertainty in the range-range rate plane, it is hard to verify if measurements belong to the confidence region, except in high-intensity cases. In the Swarm C case (not shown), only high-intensity manoeuvres are detected, with varying success for medium-intensity manoeuvres.

\subsection{Algorithm 2  simulated results}
\label{sec:six_three}

\begin{table*}[h]
	\centering
	\begin{tabular}{|l|l|l|l|l|l|l|}
		\hline
		\thead{Case / Metrics}     & \thead{P1M ($ \% $)} & \thead{P5M ($ \% $)} & \thead{P8M ($ \% $)} & \thead{P1D ($ \% $)} & \thead{P5D ($ \% $)} & \thead{P8D ($ \% $)} \\
		\hline
		\bm{$ FP $} & 10        & 1  &  0 & 3 & 0 & 0     \\
		\hline
		\bm{$ L-12h $} & 100        & 72  & 0 & 70 & 40 & 0     \\
		\hline
		\bm{$ L-6h $} & 100        & 42  & 0 & 50 & 0 & 0     \\
		\hline
		\bm{$ L-2h $} & 0        & 0  & 0 & 0 & 0 & 0     \\
		\hline
		\bm{$ M-12h $} & 100       & 96  & 0 & 100 & 84 & 0    \\
		\hline
		\bm{$ M-6h $} & 100        & 100  & 95 & 100 & 100 & 95     \\
		\hline
		\bm{$ M-2h $} & 100        & 98  & 0 & 100 & 96 & 0     \\
		\hline
		\bm{$ H-12h $} & 100        & 100  & 100 & 100 & 100 & 100    \\
		\hline
		\bm{$ H-6h $} & 100       & 100  & 100 & 100 & 100 & 100     \\
		\hline
		\bm{$ H-2h $} & 100        & 100  & 100 & 100 & 100 & 100     \\
		\hline
	\end{tabular}
	\caption{Sentinel-1A reachability analysis with Algorithm 2 and probability from optimal control distance metrics. FP=false positives, L=low, M=medium, H=high. Tangential case \label{table:three}}
\end{table*}

The results are presented in Table \ref{table:three}. The metrics computed from Algorithm 2 detect all high- and (except P8M and P8D) medium-intensity manoeuvres. P1M is the most sensitive algorithm being able to detect even some low-intensity manoeuvres. However, it has a non-negligible rate of false positives (a false positive is defined as a non-manoeuvred case from the Monte Carlo simulation being detected with $ p\geq\ 50\% $). P1D is only slightly less sensitive and reduces the number of false positives. Other metrics seen to perform worse than P1M and P1D.
In the Swarm-C case (not shown), only high-intensity manoeuvres are detected, with varying success for medium-intensity ones. As in Section \ref{sec:six_one},  main causes are long gaps without measurement and having less radar measurements right after the manoeuvre.
\section{Results for real scenarios}
\label{sec_seven}
Given the algorithms already presented and validated through simulated scenarios on previous sections, this section presents the results obtained when they were tested on real-world data, for satellites of the Sentinel and Swarm family,  and also for TerraSAR-X and TanDEM-X satellites. The data used and the sources were, for the satellites' orbits, OEM data, this is, accurate position and velocity information of the satellites under study (with precisions one order of magnitude better than the radar data, i.e., with position error of about 1 meter), provided by ESA/ESOC and DLR/GSOC. For radar data, real tracks from the Spanish survey radar S3TSR were used, with the necessary uncertainty information for the algorithms. Finally, for testing purposes,  manoeuvre data, providing accelerations in a local reference frame as well as the duration, were provided by ESA/ESOC and DLR/GSOC.

First, the selected scenarios are presented and briefly described in Section \ref{sec:seven_one}. In Section \ref{sec:seven_two} the particularities of the dynamical modelling are detailed, followed by a brief comment on the data consistency check (Section \ref{sec:seven_three}). The last subsections (\ref{sec:seven_four})--(\ref{sec:seven_five}) contain the numerical results of UKF and reachability analysis using the real data.

\subsection{Real testing scenarios}
\label{sec:seven_one}

\begin{table*}
	\centering
	\begin{tabular}{|l|l|l|l|l|l|}
		\hline
		\thead{Scenario} & \thead{Epoch of most \\ intense man.} & \thead{Initial epoch} & \thead{Final epoch} & \thead{\bm{$\#$} of segments} \\
		\hline
		\bm{$ 1 $}        & 14-Aug-2019 23:11:03  &  06-Aug-2019 06:52:43 & 16-Aug-2019 18:09:54 & 12    \\
		\hline
		\bm{$ 2 $}        & 29-Jan-2020 23:11:08 &  23-Jan-2020 06:36:14 & 31-Jan-2020 07:09:05 & 8   \\
		\hline
		\bm{$ 3 $}       & 21-May-2020 01:06:08  &  17-May-2020 06:27:59 & 25-May-2020 18:01:43 & 10    \\
		\hline
		\bm{$ 4 $}        & 19-Aug-2020 22:29:28  &  19-Aug-2020 17:45:21 & 27-Aug-2020 18:18:18 & 9    \\
		\hline
		\bm{$ 5 $}        & 14-Aug-2019 23:59:53  &  05-Aug-2019 07:00:17 & 17-Aug-2019 07:00:15 & 13    \\
		\hline
		\bm{$ 6 $}       & 19-Sep-2019 00:57:21  &  10-Sep-2019 18:00:59 & 22-Sep-2019 07:00:19 & 12    \\
		\hline
		\bm{$ 7 $}        & 20-May-2020 22:37:36  &  12-May-2020 18:09:12 & 28-May-2020 17:36:30 & 17    \\
		\hline
		\bm{$ 8 $}        & 17-Jun-2020 22:05:17  &  09-Jun-2020 17:36:29 & 24-Jun-2020 18:00:59 & 18    \\
		\hline
		\bm{$ 9 $}        & 05-Feb-2020 16:09:00  &  03-Feb-2020 22:25:25 & 13-Feb-2020 22:25:22 & 8    \\
		\hline
		\bm{$ 10 $}       & 10-Sep-2020 17:10:27  & 03-Sep-2020 11:22:01 & 10-Sep-2020 22:25:35 & 6    \\
		\hline
		\bm{$ 11 $}        & 18-Sep-2019 16:59:55  & 11-Sep-2019 11:01:44 & 21-Sep-2019 22:25:30 & 5    \\
		\hline
		\bm{$ 12 $}        & 17-Sep-2020 16:10:04  & 08-Sep-2020 11:22:02 & 19-Sep-2020 10:51:42 & 9    \\
		\hline
		\bm{$ 13 $}       & 15-Jul-2020 17:02:15  & 06-Jul-2020 01:47:22 & 20-Jul-2020 00:05:35 & 12    \\
		\hline
		\bm{$ 14 $}        & 22-Aug-2020 00:24:06  & 17-Aug-2020 17:49:25 & 24-Aug-2020 06:48:33 & 4    \\
		\hline
		\bm{$ 15 $}       & 22-Aug-2020 00:24:07  & 16-Aug-2020 18:06:50 & 22-Aug-2020 17:58:07 & 4    \\
		\hline
	\end{tabular}
	\caption{List of real testing scenarios. A segment is defined as the elapsed time between a radar track and the next.}
	\label{table:four}
\end{table*}

The list of scenarios is in Table \ref{table:four} with the corresponding epochs. The satellites used to create these scenarios are: Sentinel-1A, Sentinel-1B, Sentinel-2A, Sentinel-2B, Swarm-C, TanDEM-X and TerraSAR-X, but the relation between satellite and scenario has been omitted here. Scenarios are divided in a number of segments, which start and end at consecutive radar tracks. For the purposes of testing the RA algorithms, these segments are considered and processed individually (using the precise orbits to determine the starting point for each segment), whereas the filter runs for a full scenario, processing each segment consecutively.

\begin{figure*}[h]
	\centering
	\includegraphics[scale=0.8]{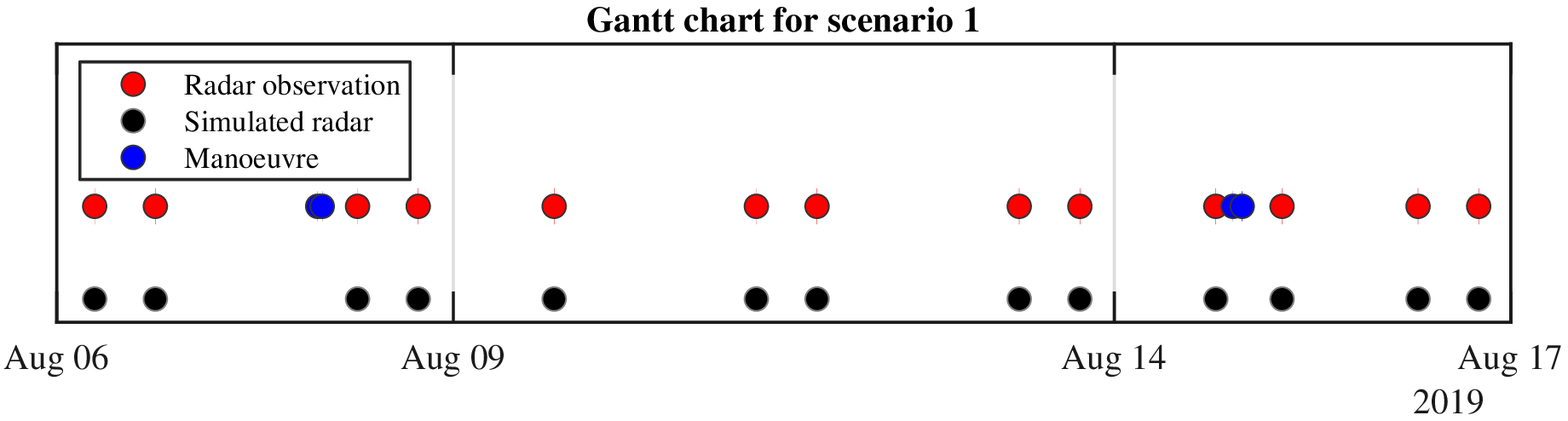}
	\caption{Gantt chart for one real scenario. This is representative of the frequency of radar tracks for the orbits under study.}
	\label{fig:gant1}
\end{figure*}

A Gantt-like representation was produced to exemplify how the scenarios are generally distributed, see Figure \ref{fig:gant1} for Scenario 1, with the radar observation (red) and manoeuvres (blue). For completeness, the plots include the simulated radar observation (black circles), which in some cases reveal missing radar tracks from the real data.

The simulations studied in Section \ref{sec:six} were carried out for a single type manoeuvres only, namely a
uniform acceleration in a single time segment, either IP, OOP or hybrid. However, real operation of satellites show that orbit corrections are generally a combination of more than one type, in consecutive and close time
segments (usually low impulse IP and medium impulse OOP). These consecutive manoeuvres may have a bigger impact on the orbit by increasing its detectability. 

\subsection{Dynamical modelling}
\label{sec:seven_two}

Since real orbits are subject to multiple complex perturbations, an important initial step is to determine the dynamical model for each satellite. In order to use the dataset provided for the real testing scenarios several improvements had to be made to the modelling of the dynamics with respect to the one used in Section \ref{sec:six}. The most important ones are the changes of the Earth gravity and atmosphere models. For one, the degree and order of the harmonics for the earth gravity field has been considerably increased. Testing has been done to discern the relation between the computational costs of increasing this parameter against the changes in the simulation error (measured with respect OEM data). 

\begin{figure}[ht!]
	\centering
	\includegraphics[scale=0.68]{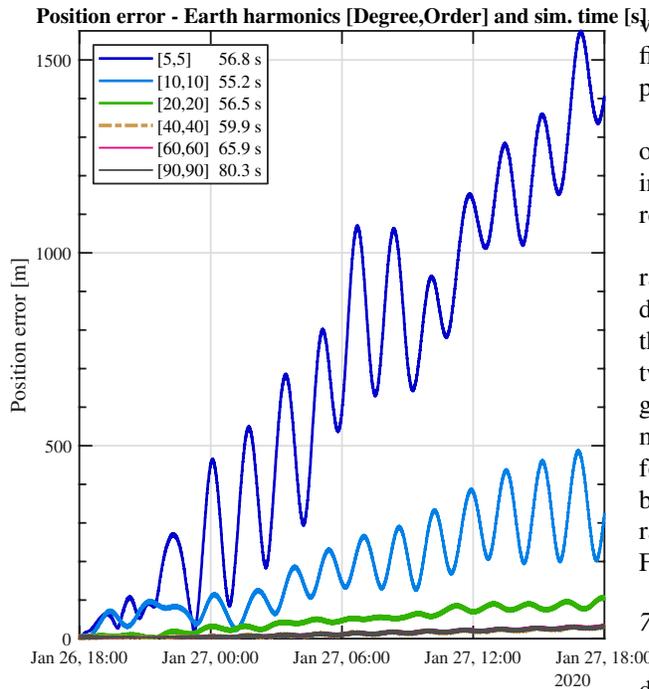}
	\caption{Position error ($ m $) for Sentinel-1A in a 24-hour segment, comparing the effect of increasing the Earth
		harmonics.}
	\label{fig:errorGravModel}
\end{figure}

Given the uncertainty of the data, values of the harmonic's degree/order above 40 have a negligible effect and can be discarded, as it would significantly slow down the computation without any relevant benefit. To justify this, a comparison has been made with increasing degree/order of the harmonics, see Figure \ref{fig:errorGravModel} for a representation of the position error evolution (against OEM data points) along a 24-hour simulation for Sentinel-1A. The time required to simulate is in the legend, where up to degree/order 40 it is affected very little when considering the great reduction in prediction error. Going above this value has a measurable effect in the time required, but with almost no impact on the error (there is even some small random increase possibly due to other perturbations and misfits). This result supports the decision to keep the harmonics only up to degree/order 40. Just to make sure that these results hold for lower orbits, a similar test has been done for a 1-day simulation interval with the OEM data of Swarm-C, which despite not being shown here points to the same conclusions.

\subsubsection{Atmosphere model and other orbital perturbations}

The atmosphere model has been changed from a Modified Harris Priester model (static atmosphere) on the simulated scenarios to the 2001 Naval Research Laboratory Mass Spectrometer and Incoherent Scatter Radar Exosphere of the MSIS model \citep{picone2002nrlmsise}, or NRLMSISE-00 model (also used by ESA for prediction and orbit determination). This last model is empirical and needs real weather data to compute the neutral atmosphere from the surface to the lower exosphere. The model feeds from the MSAFE (Marshal Solar Activity Future Estimation) data implemented inside OREKIT, which provides the mean and daily solar flux and geomagnetic indices. From the carried out testing (results omitted here for brevity), the MSAFE atmosphere provides a much better error estimation (against OEM data), but is considerably more expensive to compute, as it requires to perform interpolations from the empirical data to compute density values. These tests show us that, once the Earth shape has been modelled accurately (degree/order of 40), the gain from a more advanced atmosphere model is also very significant,  even close to one order of magnitude. Although it is not shown here, another comparison has been made to measure the relative importance of the solar radiation pressure in simulations of these periods (the order of days) when the other forces are modelled with as much precision as possible, to conclude that this perturbation's relative significance is minimal (the effect is masked by other modelling errors for the length of these simulations).

\subsubsection{Satellite parameters}

The last consideration has been the model of the satellite itself, which is defined by the drag coefficient ($ C_D $), the frontal area for the drag force ($ S $) and the area affected by the solar radiation pressure ($ S_{SRP} $). The mass of the satellite is considered different for each case and has been chosen to be the mean value between the wet and dry mass of each satellite. Simple as the model is, the 3 parameters are adjusted to get a good fit with the real data provided. In order to do this, an iterative optimization process has been performed, with enough iterations so that the changes at the end are sufficiently small to consider that it has reached a minimum (normally at around 4-5 iterations). The algorithm is quite expensive computationally speaking, as for each of the optimization steps the orbit must be simulated several times for it to find a solution. This method has another drawback, and it is that depending on the chosen segment the final values may vary slightly. Although marginal, this effect can be palliated if the final values are averaged between different segments.

\subsection{Data consistency check}
\label{sec:seven_three}

Several sanity checks were performed to ensure the consistency of the different data sources and with the propagators, namely, verifying: that the precise orbit replicate, approximately, the radar measurements (to measurement and orbit error); that the propagators do not have much error with respect to the precise orbits and the measurement in the absence of manoeuvres; that the manoeuvre file was consistent with the OEMs, which can be verified by the error of the propagators growing rapidly in the presence of manoeuvres.

These checks are essential to ensure that false positives or false negatives are not in fact detecting inconsistencies in the data sets, and although very extensive, only general results will be mentioned here. 

\begin{figure}[h!]
	\centering
	\includegraphics[scale=0.95]{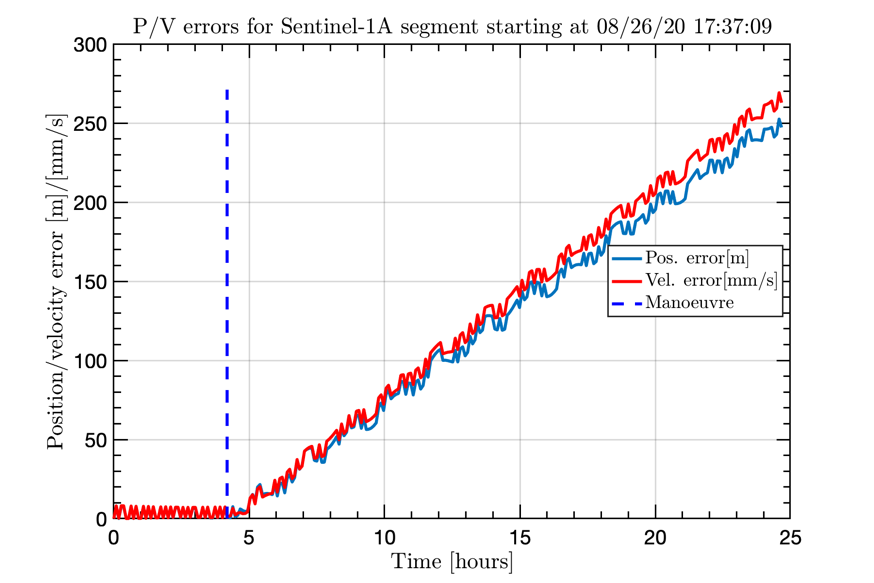}
	\includegraphics[scale=0.95]{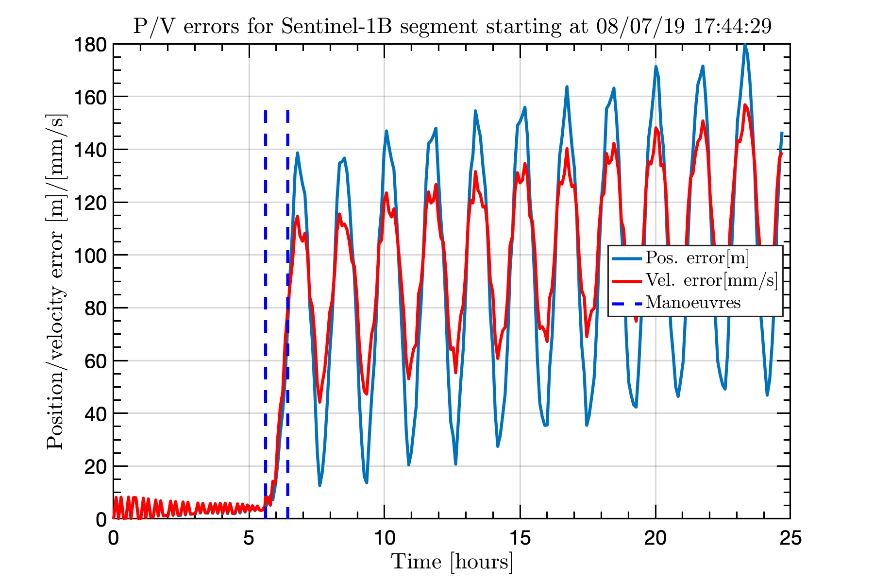}
	\caption{Examples of how manoeuvres make prediction errors (with respect to OEMs) grow large.}
	\label{fig:manError}
\end{figure}

The first check (consistency between OEM files and radar data) shows that the range differences are in the order of the combined error of the radar measurements and the OEMs themselves. The second check, consistency between our propagator and the OEM data, shows that in general the propagation performs well (errors about 60 meters maximum for a 24-hour period of propagation, as found when adjusting the values of the dynamics model), but as expected in the presence of manoeuvres errors grow rapidly, so that the last check is confirmed as well (see Figure \ref{fig:manError} for an example).

\subsection{UKF real results}
\label{sec:seven_four}

This section tests the developed UKF against real radar data of manoeuvring satellites. The filter was improved, including also the ballistic coefficient  in the estimation. In addition, since the metric $\Psi$ has a close relation to the Mahalonobis distance, the formula (\ref{eq:prmd}) was directly used to derive a manoeuvre probability.

\begin{figure}[!t]
	\centering
	\includegraphics[scale=0.95]{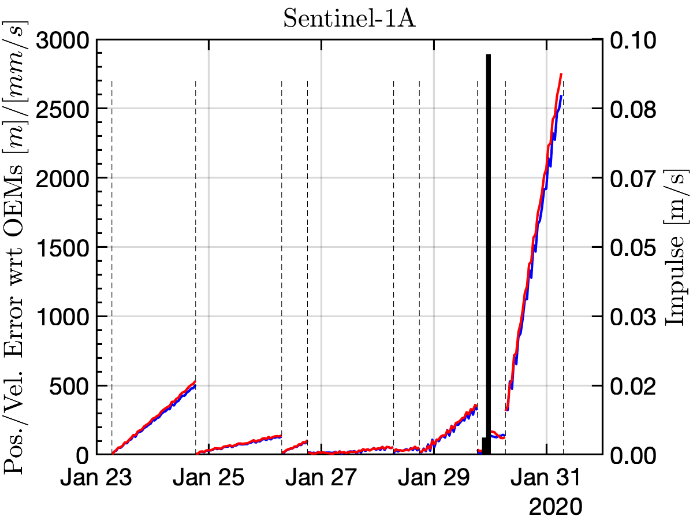}
	\includegraphics[scale=0.95]{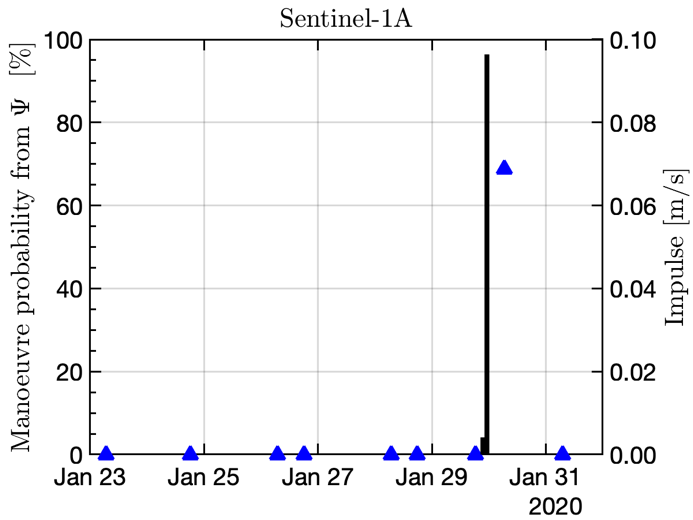}
	\includegraphics[scale=0.95]{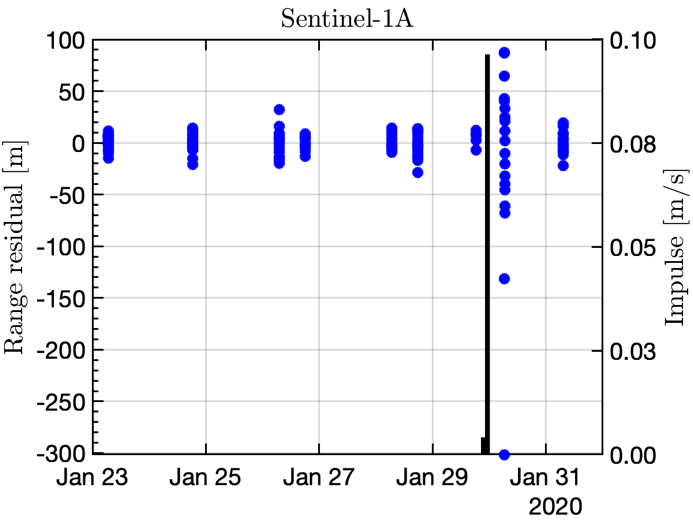}
	\caption{UKF results for one scenario of Sentinel-1A}
	\label{fig:ukfResSc2}
\end{figure}

The smoothed prediction errors with respect to OEMs, the smoothed range residuals as well as the manoeuvre detection probability derived from $ \Psi $ are shown for one scenario of Sentinel-1A in Figure \ref{fig:ukfResSc2}. The filter error increases slightly at the beginning and then the filter converges; later, after the manoeuvre, errors start to increase considerably. In any case, the steep increment in the residuals allows the detection of this manoeuvre by the filter metric.

Space restrictions do not allow for a more comprehensive revision of the results, but the following conclusions were derived. The estimation of the ballistic coefficient does not vary much and does not seem to have significant impact. It was observed that the filter was, in general, well-behaved even in the presence of manoeuvres. When there are no manoeuvres, it tends to converge, albeit sometimes slowly, to errors of the order of just a few hundred meters compatible with the errors of the propagator (for long propagations). However, there were some instances of rapid increase of errors due to the presence of manoeuvres. The manoeuvre detection metric $ \Psi^* $ allowed for the computation of a probability but was not very sensitive. Occasionally it detected a manoeuvre after one or two radar tracks have passed. The need of a combined manoeuvre detection tool integrated in the filter became clear: then, the filter can react to a manoeuvre by increasing the process noise (covariance inflation) and thus take into account the presence of the unknown manoeuvre. In addition, it would allow for longer smoothings reducing the risk of using a segment with a manoeuvre for a long smoothing. This filter is shown in Section~\ref{sec:MDF}.

\subsection{Algorithm 1 real results}
\label{sec:seven_five}

For each segment of the real scenarios, Algorithm 1 of Section \ref{sec:cinco_uno} was implemented. Initialization of the predicted orbit is done using the precise orbit information at the closest OEM point right after a radar track, and it is stopped at the attributable epoch of the next track (the middle of it). Each of these simulations are done using Taylor differential algebra and an assumed covariance matrix for the initialization information, so that the uncertainty of the prediction is known. This, combined with the confidence of the radar attributable can be used to compute a cloud of points and from that, the $\mathrm{PR}_{\mathrm{MD}}$ metric  defined  in (\ref{eq:prmd}); here MD is calculated from range and range-rate only.

The results are summarized and analysed, with each satellite being grouped in Table \ref{table:alg1real}. The results of group $\textit{1-4}$, are quite positive, as there is only one false negative and no false positives, out of 33 cases. Thus, out of 6 manoeuvres, 5 were detected and with high probability in general. This is a rather satisfactory result, as the only manoeuvre that was not detected is also the smallest. Group $\textit{5-8}$ cases contain 5 manoeuvres, but only 1 is detected (another one results in a non-zero probability but with low confidence). There are also 3 false negatives, from 55 total cases. Group $\textit{9-10}$ results are not very reassuring as no manoeuvres are detected; and with one false positive (from 22 cases). The 2 manoeuvres present in the $\textit{11-12}$ scenario are detected, with no false negatives (from 13 cases). Group $\textit{13}$ presents no detected manoeuvres (out of one) but no false negatives. For group $\textit{14-15}$ (which is a high-manoeuvring case) only $ 25 \% $ of manoeuvres are detected.

In global, the results are in need of improvement, as only about $ 40 \% $ of manoeuvres are detected, even with some of them being rather intense. The rate of false positives is quite good on the other hand. Analysing the results, the causes of errors were as follows; from a total of 5 false positives, all except one present less than 10 plots in the track following the manoeuvre. Thus, the main cause of false positives is tracks with fewer plots. From a total of 14 false negatives, all except one were segments of length equal to or larger than one day. Thus, the main cause of false negatives is longer propagations accumulating additional propagation error.	 Sometimes these longer propagation periods are due to missed radar observations right after the manoeuvre. The most challenging scenarios were those with TanDEM-X/TerraSAR-X due to the abundance of manoeuvres and the scarcity of data. This algorithm performed poorly in those scenarios compared with the others.

\begin{table}
	\centering
	\begin{tabular}{|l|l|l|l|l|}
		\hline
		\thead{Scenario} & \thead{\bm{$\#$} man.}  & \thead{\bm{$\%$} man. \\ detected} & \thead{\bm{$\%$} false \\ positives} & \thead{\bm{$\%$} false \\ negatives} \\
		\hline
		1 & 2 & 100 & 0 & 0 \\
		\hline
		2 & 1 & 100 & 0 & 0 \\
		\hline
		3 & 1 & 100 & 0 & 0 \\
		\hline
		4 & 2 & 50 & 0 & 50 \\
		\hline
		\textbf{1-4} & 6 & 83.3 & 0 & 16.7 \\
		\hline
		5 & 2 & 50 & 0 & 50 \\
		\hline
		6 & 1 & 0 & 0 & 100 \\
		\hline
		7 & 1 & 0 & 6.25 & 100 \\
		\hline
		8 & 1 & 0 & 11.76 & 100 \\
		\hline
		\textbf{5-8} & 5 & 20 & 5.45 & 80 \\
		\hline
		9 & 1 & 0 & 9.09 & 100 \\
		\hline
		10 & 1 & 0 & 0 & 100 \\
		\hline
		\textbf{9-10} & 2 & 0 & 4.55 & 100 \\
		\hline
		11 & 1 & 100 & 0 & 0 \\
		\hline
		12 & 1 & 100 & 0 & 0 \\
		\hline
		\textbf{11-12} & 2 & 100 & 0 & 0 \\
		\hline
		13 & 1 & 0 & 0 & 100 \\
		\hline
		14 & 4 & 25 & N/A & 75 \\
		\hline
		15 & 4 & 25 & N/A & 75 \\
		\hline
		\textbf{14-15} & 8 & 25 & N/A & 75 \\
		\hline
		\textbf{All} & \textbf{16} & \textbf{41.66} & \textbf{2.98} & \textbf{58.34 }\\
		\hline
	\end{tabular}
	\caption{Algorithm 1 summarized results.}
	\label{table:alg1real}
\end{table}

\subsection{Algorithm 2 real results}
\label{sec:seven_six}
\begin{table*}
	\centering
	\begin{tabular}{|l|l|l|l|l|l|l|l|}
		\hline
		\thead{Scenario} & \thead{\bm{$\#$} man.}  & \thead{P1M \bm{$\%$} \\ man. \\ detected} & \thead{P1M \bm{$\%$} \\ false \\ positives} & \thead{P1M \bm{$\%$} \\ false \\ negatives}  & \thead{P1D \bm{$\%$} \\ man. \\ detected} & \thead{P1D \bm{$\%$} \\ false \\ positives} & \thead{P1D \bm{$\%$} \\ false \\ negatives} \\
		\hline
		1 & 2 & 100 & 0 & 0 & 100 & 0 & 0  \\
		\hline
		2 & 1 & 100 & 14.29 & 0 & 100 & 14.29 & 0  \\
		\hline
		3 & 1 & 100 & 0 & 0 & 100 & 0 & 0  \\
		\hline
		4 & 2 & 50 & 0 & 50 & 50 & 0 & 50  \\
		\hline
		\textbf{1-4} & 6 & 83.3 & 3.03 & 16.7 & 83.3 & 3.03 & 16.7  \\
		\hline
		5 & 2 &50 &0 &50 &50 &27.3 & 50  \\
		\hline
		6 & 1 & 0 & 0 & 100 & 0 & 9.1 & 100  \\
		\hline
		7 & 1 & 0 & 12.50 & 100 & 0 & 12.50 & 100  \\
		\hline
		8 & 1 & 0 & 11.76 & 100 & 0 & 11.76  &100  \\
		\hline
		\textbf{5-8} & 5 & 20 & 7.27 & 80 20 & 5.45 & 80  \\
		\hline
		9 & 1 & 0 & 20 & 100 & 0 & 20 & 100  \\
		\hline
		10 & 1 & 0 & 25 & 100 & 0 & 25 & 100  \\
		\hline
		\textbf{9-10} & 2 & 0 & 22.73 & 100 & 0 & 22.73 & 100  \\
		\hline
		11 & 1 & 100 & 0 & 0 & 100 & 0 & 0  \\
		\hline
		12 & 1  & 0  & 10  & 100  & 0  & 10  & 100  \\
		\hline
		\textbf{11-12} & 2 & 50 & 7.69 & 50 & 50 & 7.69 & 50  \\
		\hline
		13 & 1 & 0 & 0 & 100 & 0 & 0 & 100 \\
		\hline
		14 &4 &50 & N/A &50 &25& N/A & 75  \\
		\hline
		15 & 4 & 100 & N/A & 0 & 100 & N/A & 0 \\
		\hline
		\textbf{14-15} & 8 &75 &N/A &25 &62.5 &N/A &37.5  \\
		\hline
		\textbf{All} & 24 & 54.16& 45.84& 8.21& 50 &7.46& 50  \\
		\hline
	\end{tabular}
	\caption{Algorithm 2 summarized results.}
	\label{table:alg2real}
\end{table*}
The summarized results can be seen in Table \ref{table:alg2real} (results have to consider there are 134 segments without manoeuvres). 

The results of group $\textit{1-4}$ are quite positive for both P1M and P1D, as there is only one false negative and few false positives, out of 33 cases. Thus, out of 6 manoeuvres, 5 were detected. Only the smallest manoeuvre was not detected (similarly to the Algorithm 1 results). Group $\textit{5-8}$ cases contain 5 manoeuvres, but only 1 is detected for both P1M and P1D (another one is felt but with low confidence). There are also some false positives, from 55 total cases, with P1M obtaining a slightly higher rate of false positives. Group $\textit{9-10}$ results are quite bad for both metrics as no manoeuvres are detected; and with a rather high rate of false positive (from 22 cases). One of the 2 manoeuvres present in group $\textit{11-12}$ scenarios is detected, with some false negatives (from 13 cases). Results are the same for both metrics. Group $\textit{13}$ presents no detected manoeuvres (out of one) but no false negatives, for both metrics. Finally, group $\textit{14-15}$ gives 8 segments, all of them with manoeuvres; out of these, 5 are detected for P1D and 6 for P1M. In global, the results can be considered positive, as more than half the manoeuvres are detected with a low rate of false negatives for both metrics, but in need of improvement. P1M seemed to perform better than P1D, with a minimal increase in false negatives. Analysing the results, the causes of errors are similar as for Algorithm 1: From a total of 11 (10) false positives for P1M (resp., P1D), all except two (resp., one) present less than 10 plots in the track following the manoeuvre. Thus, the main cause of false positives is tracks with fewer plots. From a total of 11 (12) false negatives for P1M (resp. P1D), all except one were segments of length equal to or larger than one day. Thus, the main cause of false negatives is longer propagations accumulating additional	propagation error. Sometimes these longer propagation periods are due to missed radar observations right after the manoeuvre. The most challenging scenarios are those with TanDEM-X/TerraSAR-X due to the abundance of manoeuvres and the scarcity of data. As opposed to Algorithm 1, this algorithm performs excellently in those scenarios compared with the other.

\section{A manoeuvre detection filter using reachability analysis}
\label{sec:MDF}
 The results of the previous section shed light on the need of a combined manoeuvre detection tool-filter, that is able both to predict an orbit and detect and take into account manoeuvres. This section presents our work on such a filter, which we call the Manoeuvre Detection Filter (MDF).
 The idea is to follow the scheme of the UKF of \ref{sec:seven_four}, combined with the Algorithm 1 of \ref{sec:seven_five}, to detect manoeuvres.
Once a manoeuvre is detected, the idea of covariance inflation is followed; thus, the state covariance is increased up to the point where a manoeuvre is no longer detected, which would imply that the uncertainty of the state is able to include the possibility of such a manoeuvre having been performed. The reasons to choose Algorithm 1 instead of 2 are that it performs better for scenarios 1 to 13, which are the ones best suited to the filter. It is hopeless to expect the filter to perform well in scenarios with too many manoeuvres and few radar tracks. In addition, Algorithm 1 fits quite well with the philosophy of the filter: the unscented transform can be used to estimate the state covariance by using the attributable as a ``virtual measurement'' used only for purposes of manoeuvre detection, but not for updating the state. This considerably reduces the computational burden. Also Algorithm 1 gave less false positives in the real testing.

To be more precise, using the UKF notation of Section~\ref{sect-ukfalg} and skipping the unchanged steps, the MDF algorithm is:

\begin{enumerate}
	\item Start from the previous estimate of the state and the covariance of its error ($ \bm{\hat{x}}_0 $ and $ \bm{\hat{P}}_0 $ at first).
	\item Compute the attributable time $t_i^\mathcal{A}$, the values of the next track $\bm{y}_i^\mathcal{A}$ and its covariance $\bm{\Sigma}_i^\mathcal{A}$.
	\item Compute the sigma-points of the unscented transform.
	\item Propagate all the sigma points using numerical integration until the attributable time $\tilde{\bm{x}}_{i}^{\left(j\right)}$.
	\item Compute the weighted mean and the covariance matrix of the transformed sigma-points: $\bar{\bm{x}}_i$ and $\bm{\bar{\mathrm{P}}}_i$.
	\item Compute the probability of manoeuvre $ p_i $ by transforming the sigma-points (using the observation equation) to get the predicted observation $\hat{\bm{y}}_{i}$, the residuals $\bm{\nu}_i=\bm{y}_{i}^\mathcal{A} -\hat{\bm{y}}_{i}$, and the observation covariance $\bm{\mathrm{S}}_i$ to compute the $\mathrm{MD}_i$ and the associated probability $ p_i $ (using Equation \ref{eq:prmd}):
	\begin{equation*}
		\bm{\mathrm{S}}_i = \Sigma^{2n}_{j=0} w_c^{\left(j\right)} \left(\tilde{\bm{y}}_{i}^{\left(j\right)} - \hat{\bm{y}}_{i} \right)\left(\tilde{\bm{y}}_{i}^{\left(j\right)} - \hat{\bm{y}}_{i} \right)^T + \bm{\Sigma}_i^\mathcal{A}
	\end{equation*}
	\begin{equation*}
		\mathrm{MD}_i = \sqrt{\bm{\nu}_i^T\bm{\mathrm{S}}_i^{-1}\bm{\nu}_i}.
	\end{equation*}
	
	\item If $ p_i \ge 0,5 $, multiply the covariance of the predicted state $\bm{\bar{\mathrm{P}}}_i$ by 2 and return to Point 6; otherwise, continue. 
	\item Transform the sigma-points (using the observation equation) and calculate the predicted observation, the residuals, and the observation covariance.
	\item Calculate the predicted observation covariance and the residuals.
	\item Compute the Kalman gain and update the state estimate.
	\item Return to step 1 and continue propagating.
\end{enumerate}

In addition, a ``long smoothing'' is implemented: if no manoeuvre is detected, a  smoothing is performed backwards until the previous radar track, and again forwards. 

\begin{figure}[ht!]
	\centering
	\includegraphics[scale=1]{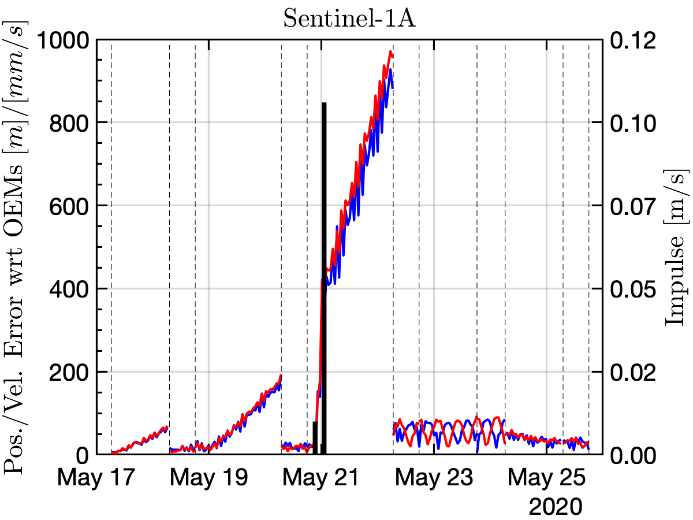}
	\includegraphics[scale=1]{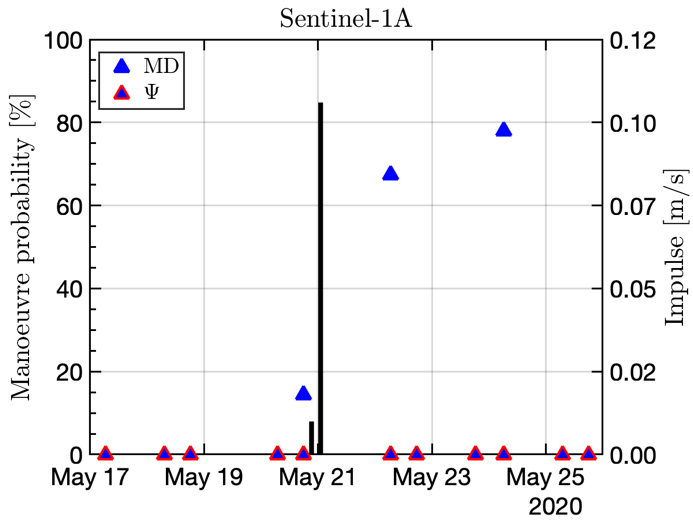}
	\includegraphics[scale=1]{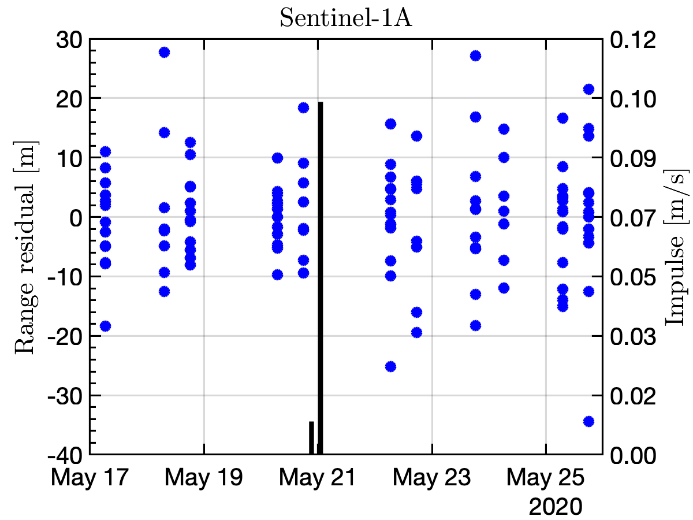}
	\caption{MDF results for one scenario of Sentinel-1A.}
	\label{fig:mdfResSc2}
\end{figure}

Figure \ref{fig:mdfResSc2} has one example of the MDF results. In this case, the filter performs initially quite well thanks to the long smoothing, and the manoeuvre is detected. Even though the initial spike after the manoeuvre is quite large, the filter recovers quite quickly thanks to the inflation mechanism; comparing with the UKF of Figure~\ref{fig:ukfResSc2} the behaviour is much improved. However, a false positive also happens before the end, but it does not impact the filter's performance. Note that the UKF metric did not detect manoeuvres in this segment.

The following conclusions can be derived from the results obtained from the MDF, which cannot be shown here due to lack of space. It can be observed that the filter is, in general, well-behaved even in the presence of manoeuvres and detects many of them. When there are no manoeuvres, it tends to converge, quicker than the UKF without long smoothing, to errors of the order of just a few hundred meters compatible with the errors of the propagator.  As in the UKF, there are some instances of rapid increase of errors due to the presence of manoeuvres, particularly when undetected. The MD metric allows for detection of many manoeuvres but also produces a considerable number of false positives. It is not as straightforward as for the RA algorithms to obtain fair statistics, since it is unclear if a detection after one or two tracks should be considered a true or false positive; this is due to the sequential nature of this algorithm, which considers scenarios as a whole, instead of processing segments separately. Thus, the history of each scenario influences the results in several ways. Covariance inflation works well in most cases but in some instances, it might be too large, inducing large errors in the state.

\section{Conclusions and future work}
\label{sec:eight}
Several methods for the detection of manoeuvres in LEO from radar data have been presented, based on UKF, attributable theory and reachability analysis. Initial simulation results showed that the filter did not detect manoeuvres unless they are rather intense, whereas the reachability approach was more sensitive at the price of longer computation times. When tested with real data, the results, while not bad, are in need of improvement. Since the quality of the data was verified, the main identified difficulty was the scarcity of measurements (low number of tracks resulting in long propagation times without information and/or low number of plots in some cases), due to the fact of having a single surveillance radar, the Spanish survey radar S3TSR, as the source of data. Sometimes, more than 24 hours or more went without a measurement and the manoeuvres were not very intense. Future ideas to address this challenge include improvements in the propagator, in the description of uncertainty (e.g. the use of Gaussian mixtures to better describe the error distribution after long propagations), and the use of the intensive surveillance mode of the radar, which can provide hundreds of plots for a single track.

From an innovation point of view, the results can be considered of interest, since most of the results in the literature depend on having numerous measurements (oftentimes, almost continuous data is assumed, which is only realistic in GEO with optical sensors), which was not the case here. The final aim  is to have these algorithms integrated in the S3T Cataloguing System in order to provide routine automatic manoeuvre detection capabilities, but all methods presented here can be refined and extended in many directions. Beyond obvious extensions or refinement of the algorithms, there is much to be gained from additional measurements (additional stations, more regular radar tracks, SLR measurements, etc.) as it was shown that RA algorithms under-perform when the number of radar plots are lower. It would be of great interest to draw a set of minimum measurement requirements for the successful application of the manoeuvre detection algorithms.

\section*{Acknowledgements}
The authors would like to acknowledge Centro de Desarrollo Tecnol\'ogico e Industrial (CDTI) for providing funding for this work under the S3T program and in particular  Cristina P\'erez Hern\'andez for her help and support. The authors are also grateful to Ralph Kahle from DLR for providing the orbital and manoeuvre data for TerraSAR-X and TanDEM-X satellites. Finally, Universidad de Sevilla is acknowledged for partial funding under V PPIT-US and VI PPIT-US.
\bibliographystyle{model5-names}
\biboptions{authoryear}
\bibliography{refs}

\begin{thebibliography}{33}
\expandafter\ifx\csname natexlab\endcsname\relax\def\natexlab#1{#1}\fi
\providecommand{\url}[1]{\texttt{#1}}
\providecommand{\href}[2]{#2}
\providecommand{\path}[1]{#1}
\providecommand{\DOIprefix}{doi:}
\providecommand{\ArXivprefix}{arXiv:}
\providecommand{\URLprefix}{URL: }
\providecommand{\Pubmedprefix}{pmid:}
\providecommand{\doi}[1]{\href{http://dx.doi.org/#1}{\path{#1}}}
\providecommand{\Pubmed}[1]{\href{pmid:#1}{\path{#1}}}
\providecommand{\bibinfo}[2]{#2}
\ifx\xfnm\relax \def\xfnm[#1]{\unskip,\space#1}\fi
\bibitem[{Andersson et~al.(2019)Andersson, Gillis, Horn
  et~al.}]{andersson2019casadi}
\bibinfo{author}{Andersson, J.~A.}, \bibinfo{author}{Gillis, J.},
  \bibinfo{author}{Horn, G.} et~al. (\bibinfo{year}{2019}).
\newblock \bibinfo{title}{Cas{A}{D}i: a software framework for nonlinear
  optimization and optimal control}.
\newblock {\it \bibinfo{journal}{Mathematical Programming Computation}\/},
  {\it \bibinfo{volume}{11}\/}\bibinfo{issue}{(1)}, \bibinfo{pages}{1--36}.
\bibitem[{Andrea and Maisonobe(2016)}]{andrea2016automatic}
\bibinfo{author}{Andrea, A.}, and \bibinfo{author}{Maisonobe, L.}
  (\bibinfo{year}{2016}).
\newblock \bibinfo{title}{Automatic differentiation for propagation of orbit
  uncertainties on {O}rekit}.
\newblock In {\it \bibinfo{booktitle}{Stardust Conference 2016}\/}.
\newblock \URLprefix
  \url{https://www.orekit.org/doc/Antolino-2016-automatic-diff-for-prop-of-orbit-uncertainties.pdf}.
\bibitem[{Armellin et~al.(2010)Armellin, Di~Lizia, Bernelli-Zazzera
  et~al.}]{armellin2010asteroid}
\bibinfo{author}{Armellin, R.}, \bibinfo{author}{Di~Lizia, P.},
  \bibinfo{author}{Bernelli-Zazzera, F.} et~al. (\bibinfo{year}{2010}).
\newblock \bibinfo{title}{Asteroid close encounters characterization using
  differential algebra: the case of {A}pophis}.
\newblock {\it \bibinfo{journal}{Celestial Mechanics and Dynamical
  Astronomy}\/},  {\it \bibinfo{volume}{107}\/}\bibinfo{issue}{(4)},
  \bibinfo{pages}{451--470}.
\bibitem[{Carpenter and D'Souza(2018)}]{carpenter2018navigation}
\bibinfo{author}{Carpenter, J.~R.}, and \bibinfo{author}{D'Souza, C.~N.}
  (\bibinfo{year}{2018}).
\newblock {\it \bibinfo{title}{Navigation Filter Best Practices}\/}.
\newblock \bibinfo{type}{Technical Report} \bibinfo{number}{TP-2018-219822}
  NASA.
\newblock \URLprefix
  \url{https://ntrs.nasa.gov/api/citations/20180003657/downloads/20180003657.pdf}.
\bibitem[{Goff(2015)}]{goff2015orbit}
\bibinfo{author}{Goff, G.~M.} (\bibinfo{year}{2015}).
\newblock {\it \bibinfo{title}{Orbit Estimation of Non-Cooperative Maneuvering
  Spacecraft}\/}.
\newblock Ph.D. thesis Air Force Institute of Technology.
\newblock \URLprefix \url{https://scholar.afit.edu/etd/202/}.
\bibitem[{Goff et~al.(2015{\natexlab{a}})Goff, Black and
  Beck}]{goff2015tracking}
\bibinfo{author}{Goff, G.~M.}, \bibinfo{author}{Black, J.~T.}, and
  \bibinfo{author}{Beck, J.~A.} (\bibinfo{year}{2015}{\natexlab{a}}).
\newblock \bibinfo{title}{Tracking maneuvering spacecraft with filter-through
  approaches using interacting multiple models}.
\newblock {\it \bibinfo{journal}{Acta Astronautica}\/},  {\it
  \bibinfo{volume}{114}\/}, \bibinfo{pages}{152--163}.
\bibitem[{Goff et~al.(2015{\natexlab{b}})Goff, Showalter, Black
  et~al.}]{goff2015parameter}
\bibinfo{author}{Goff, G.~M.}, \bibinfo{author}{Showalter, D.},
  \bibinfo{author}{Black, J.~T.} et~al. (\bibinfo{year}{2015}{\natexlab{b}}).
\newblock \bibinfo{title}{Parameter requirements for noncooperative satellite
  maneuver reconstruction using adaptive filters}.
\newblock {\it \bibinfo{journal}{Journal of Guidance, Control, and
  Dynamics}\/},  {\it \bibinfo{volume}{38}\/}\bibinfo{issue}{(3)},
  \bibinfo{pages}{361--374}.
\bibitem[{Gomez et~al.(2019)Gomez, Salmer{\'o}n, Besso
  et~al.}]{gomez2019initial}
\bibinfo{author}{Gomez, R.}, \bibinfo{author}{Salmer{\'o}n, J. M.-V.},
  \bibinfo{author}{Besso, P.} et~al. (\bibinfo{year}{2019}).
\newblock \bibinfo{title}{Initial operations of the breakthrough {S}panish
  {S}pace {S}urveillance and {T}racking {R}adar ({S}{3}{T}{S}{R}) in the
  {E}uropean context}.
\newblock In {\it \bibinfo{booktitle}{1st ESA NEO and Debris Detection
  Conference. Presented paper. Darmstadt, Germany}\/}.
\newblock \URLprefix
  \url{https://conference.sdo.esoc.esa.int/proceedings/neosst1/paper/479}.
\bibitem[{Guang et~al.(2018)Guang, Xingzi, Hanyu et~al.}]{guang2018non}
\bibinfo{author}{Guang, Z.}, \bibinfo{author}{Xingzi, B.},
  \bibinfo{author}{Hanyu, Z.} et~al. (\bibinfo{year}{2018}).
\newblock \bibinfo{title}{Non-cooperative maneuvering spacecraft tracking via a
  variable structure estimator}.
\newblock {\it \bibinfo{journal}{Aerospace Science and Technology}\/},  {\it
  \bibinfo{volume}{79}\/}, \bibinfo{pages}{352--363}.
\bibitem[{Hall and Singla(2019)}]{hall2019probabilistic}
\bibinfo{author}{Hall, Z.}, and \bibinfo{author}{Singla, P.}
  (\bibinfo{year}{2019}).
\newblock \bibinfo{title}{A probabilistic approach for reachability set
  computation for efficient space situational awareness}.
\newblock In {\it \bibinfo{booktitle}{29th AAS/AIAA Space Flight Mechanics
  Meeting, 2019}\/} (pp. \bibinfo{pages}{3001--3020}).
\newblock \bibinfo{organization}{Univelt Inc.}
\bibitem[{Holzinger and Scheeres(2009)}]{holzinger2009reachability}
\bibinfo{author}{Holzinger, M.}, and \bibinfo{author}{Scheeres, D.}
  (\bibinfo{year}{2009}).
\newblock {\it \bibinfo{title}{Reachability Analysis Applied to Space
  Situational Awareness}\/}.
\newblock \bibinfo{type}{Technical Report} Colorado University at Boulder,
  departament of aerospace engineering sciences.
\newblock \URLprefix \url{https://apps.dtic.mil/sti/pdfs/ADA531761.pdf}.
\bibitem[{Holzinger et~al.(2012)Holzinger, Scheeres and
  Alfriend}]{holzinger2012object}
\bibinfo{author}{Holzinger, M.~J.}, \bibinfo{author}{Scheeres, D.~J.}, and
  \bibinfo{author}{Alfriend, K.~T.} (\bibinfo{year}{2012}).
\newblock \bibinfo{title}{Object correlation, maneuver detection, and
  characterization using control distance metrics}.
\newblock {\it \bibinfo{journal}{Journal of Guidance, Control, and
  Dynamics}\/},  {\it \bibinfo{volume}{35}\/}\bibinfo{issue}{(4)},
  \bibinfo{pages}{1312--1325}.
\bibitem[{Jain et~al.(2019)Jain, Gu{\'e}ho, Singla et~al.}]{jain2019stochastic}
\bibinfo{author}{Jain, A.}, \bibinfo{author}{Gu{\'e}ho, D.},
  \bibinfo{author}{Singla, P.} et~al. (\bibinfo{year}{2019}).
\newblock \bibinfo{title}{Stochastic reachability analysis for the hypersonic
  re-entry problem}.
\newblock In {\it \bibinfo{booktitle}{29th AAS/AIAA Space Flight Mechanics
  Meeting, 2019}\/} (pp. \bibinfo{pages}{2455--2476}).
\newblock \bibinfo{organization}{Univelt Inc.}
\bibitem[{Julier and Uhlmann(1997)}]{julier1997new}
\bibinfo{author}{Julier, S.~J.}, and \bibinfo{author}{Uhlmann, J.~K.}
  (\bibinfo{year}{1997}).
\newblock \bibinfo{title}{New extension of the {K}alman filter to nonlinear
  systems}.
\newblock In {\it \bibinfo{booktitle}{Signal processing, sensor fusion, and
  target recognition VI}\/} (pp. \bibinfo{pages}{182--193}).
\newblock \bibinfo{organization}{International Society for Optics and
  Photonics} volume \bibinfo{volume}{3068}.
\bibitem[{Kurzhanski and Varaiya(2000)}]{kurzhanski2000ellipsoidal}
\bibinfo{author}{Kurzhanski, A.~B.}, and \bibinfo{author}{Varaiya, P.}
  (\bibinfo{year}{2000}).
\newblock \bibinfo{title}{Ellipsoidal techniques for reachability analysis:
  internal approximation}.
\newblock {\it \bibinfo{journal}{Systems \& control letters}\/},  {\it
  \bibinfo{volume}{41}\/}\bibinfo{issue}{(3)}, \bibinfo{pages}{201--211}.
\bibitem[{Li and Jilkov(2005)}]{li2005survey}
\bibinfo{author}{Li, X.~R.}, and \bibinfo{author}{Jilkov, V.~P.}
  (\bibinfo{year}{2005}).
\newblock \bibinfo{title}{Survey of maneuvering target tracking. part v.
  multiple-model methods}.
\newblock {\it \bibinfo{journal}{IEEE Transactions on Aerospace and Electronic
  Systems}\/},  {\it \bibinfo{volume}{41}\/}\bibinfo{issue}{(4)},
  \bibinfo{pages}{1255--1321}.
\bibitem[{Mainsonobe et~al.(2012)Mainsonobe, Cefola, Frouvelle et~al.}]{orekit}
\bibinfo{author}{Mainsonobe, L.}, \bibinfo{author}{Cefola, P.},
  \bibinfo{author}{Frouvelle, N.} et~al. (\bibinfo{year}{2012}).
\newblock \bibinfo{title}{Open governance of the {O}rekit space flight dynamics
  library}.
\newblock In {\it \bibinfo{booktitle}{Proceedings of the International
  Conference on Astrodynamic Tools and Techniques (ICATT)}\/} (pp.
  \bibinfo{pages}{327--343}).
\newblock \bibinfo{organization}{ESA/ESTEC} volume~\bibinfo{volume}{29}.
\bibitem[{Pastor et~al.(2020)Pastor, Escribano and
  Escobar}]{pastor2020satellite}
\bibinfo{author}{Pastor, A.}, \bibinfo{author}{Escribano, G.}, and
  \bibinfo{author}{Escobar, D.} (\bibinfo{year}{2020}).
\newblock \bibinfo{title}{Satellite maneuver detection with optical survey
  observations}.
\newblock In {\it \bibinfo{booktitle}{Advanced Maui Optical and Space
  Surveillance Technologies Conference (AMOS)}\/}.
\newblock \URLprefix
  \url{https://amostech.com/TechnicalPapers/2020/Astrodynamics/Pastor.pdf}.
\bibitem[{P{\'e}rez et~al.(2013)P{\'e}rez, Masdemont~Soler and
  G{\'o}mez~Muntan{\'e}}]{perez2013jet}
\bibinfo{author}{P{\'e}rez, D.}, \bibinfo{author}{Masdemont~Soler, J.}, and
  \bibinfo{author}{G{\'o}mez~Muntan{\'e}, G.} (\bibinfo{year}{2013}).
\newblock \bibinfo{title}{Jet {T}ransport propagation of uncertainties for
  orbits around the earth}.
\newblock In {\it \bibinfo{booktitle}{64rd International Astronautical
  Congress}\/} (pp. \bibinfo{pages}{1--8}).
\newblock \URLprefix \url{https://upcommons.upc.edu/handle/2117/21808}.
\bibitem[{Picone et~al.(2002)Picone, Hedin, Drob et~al.}]{picone2002nrlmsise}
\bibinfo{author}{Picone, J.}, \bibinfo{author}{Hedin, A.},
  \bibinfo{author}{Drob, D.~P.} et~al. (\bibinfo{year}{2002}).
\newblock \bibinfo{title}{{N}{R}{L}{M}{S}{I}{S}{E}-00 empirical model of the
  atmosphere: Statistical comparisons and scientific issues}.
\newblock {\it \bibinfo{journal}{Journal of Geophysical Research: Space
  Physics}\/},  {\it \bibinfo{volume}{107}\/}\bibinfo{issue}{(A12)},
  \bibinfo{pages}{SIA--15}.
\bibitem[{Poore et~al.(2016)Poore, Aristoff, Horwood
  et~al.}]{poore2016covariance}
\bibinfo{author}{Poore, A.~B.}, \bibinfo{author}{Aristoff, J.~M.},
  \bibinfo{author}{Horwood, J.~T.} et~al. (\bibinfo{year}{2016}).
\newblock {\it \bibinfo{title}{Covariance and Uncertainty Realism in Space
  Surveillance and Tracking}\/}.
\newblock \bibinfo{type}{Technical Report} Numerica Corporation.
\newblock \URLprefix \url{https://apps.dtic.mil/sti/pdfs/AD1020892.pdf}.
\bibitem[{Reihs et~al.(2021)Reihs, Vananti, Schildknecht
  et~al.}]{reihs2021application}
\bibinfo{author}{Reihs, B.}, \bibinfo{author}{Vananti, A.},
  \bibinfo{author}{Schildknecht, T.} et~al. (\bibinfo{year}{2021}).
\newblock \bibinfo{title}{Application of attributables to the correlation of
  surveillance radar measurements}.
\newblock {\it \bibinfo{journal}{Acta Astronautica}\/},  {\it
  \bibinfo{volume}{182}\/}, \bibinfo{pages}{399--415}.
\bibitem[{Sanchez et~al.(2019)Sanchez, Louembet, Gavilan
  et~al.}]{sanchez2019event}
\bibinfo{author}{Sanchez, J.~C.}, \bibinfo{author}{Louembet, C.},
  \bibinfo{author}{Gavilan, F.} et~al. (\bibinfo{year}{2019}).
\newblock \bibinfo{title}{An {E}vent-{T}riggered predictive controller for
  spacecraft rendezvous hovering phases}.
\newblock {\it \bibinfo{journal}{IFAC-PapersOnLine}\/},  {\it
  \bibinfo{volume}{52}\/}\bibinfo{issue}{(12)}, \bibinfo{pages}{97--102}.
\bibitem[{Schutz et~al.(2004)Schutz, Tapley and Born}]{schutz2004statistical}
\bibinfo{author}{Schutz, B.}, \bibinfo{author}{Tapley, B.}, and
  \bibinfo{author}{Born, G.~H.} (\bibinfo{year}{2004}).
\newblock {\it \bibinfo{title}{Statistical orbit determination}\/}.
\newblock \bibinfo{publisher}{Elsevier}.
\bibitem[{Siminski et~al.(2017)Siminski, Flohrer and
  Schildknecht}]{siminski2017assessment}
\bibinfo{author}{Siminski, J.}, \bibinfo{author}{Flohrer, T.}, and
  \bibinfo{author}{Schildknecht, T.} (\bibinfo{year}{2017}).
\newblock \bibinfo{title}{Assessment of post-maneuver observation correlation
  using short-arc tracklets}.
\newblock {\it \bibinfo{journal}{Journal of the British Interplanetary
  Society}\/},  {\it \bibinfo{volume}{70}\/}, \bibinfo{pages}{63--68}.
\bibitem[{Singh et~al.(2012)Singh, Horwood and Poore}]{singh2012space}
\bibinfo{author}{Singh, N.}, \bibinfo{author}{Horwood, J.~T.}, and
  \bibinfo{author}{Poore, A.~B.} (\bibinfo{year}{2012}).
\newblock \bibinfo{title}{Space object maneuver detection via a joint optimal
  control and multiple hypothesis tracking approach}.
\newblock In {\it \bibinfo{booktitle}{Proceedings of the 22nd AAS/AIAA Space
  Flight Mechanics Meeting}\/} (pp. \bibinfo{pages}{2012--159}).
\newblock \bibinfo{organization}{Univelt San Diego, CA} volume
  \bibinfo{volume}{143}.
\bibitem[{Vananti et~al.(2017)Vananti, Schildknecht, Siminski
  et~al.}]{vananti2017tracklet}
\bibinfo{author}{Vananti, A.}, \bibinfo{author}{Schildknecht, T.},
  \bibinfo{author}{Siminski, J.} et~al. (\bibinfo{year}{2017}).
\newblock \bibinfo{title}{Tracklet-tracklet correlation method for radar and
  angle observations}.
\newblock In {\it \bibinfo{booktitle}{Proc. 7th European Conference on Space
  Debris, Darmstadt, Germany}\/} (pp. \bibinfo{pages}{18--21}).
\newblock \URLprefix
  \url{https://conference.sdo.esoc.esa.int/proceedings/sdc7/paper/925}.
\bibitem[{Vazquez et~al.(2021{\natexlab{a}})Vazquez, Sanchez, Montilla
  et~al.}]{vazdebris2021}
\bibinfo{author}{Vazquez, R.}, \bibinfo{author}{Sanchez, J.},
  \bibinfo{author}{Montilla, J.~M.} et~al.
  (\bibinfo{year}{2021}{\natexlab{a}}).
\newblock \bibinfo{title}{Manoeuvre detection for near-orbiting objects}.
\newblock In {\it \bibinfo{booktitle}{Proc. 8th European Conference on Space
  Debris, Darmstadt, Germany}\/}.
\newblock \URLprefix
  \url{https://conference.sdo.esoc.esa.int/proceedings/sdc8/paper/238/SDC8-paper238.pdf}.
\bibitem[{Vazquez et~al.(2021{\natexlab{b}})Vazquez, Sanchez, Montilla
  et~al.}]{vazstardust2021}
\bibinfo{author}{Vazquez, R.}, \bibinfo{author}{Sanchez, J.},
  \bibinfo{author}{Montilla, J.~M.} et~al.
  (\bibinfo{year}{2021}{\natexlab{b}}).
\newblock \bibinfo{title}{Two manoeuvre detection probability metrics based on
  radar measurements and validated with {S}{3}{T}{S}{R} data}.
\newblock In {\it \bibinfo{booktitle}{Stardust-R 2nd Global Virtual Workshop
  (GVW-II) Book of Abstracts}\/} (pp. \bibinfo{pages}{153--156}).
\bibitem[{Wan and Van Der~Merwe(2000)}]{wan2000unscented}
\bibinfo{author}{Wan, E.~A.}, and \bibinfo{author}{Van Der~Merwe, R.}
  (\bibinfo{year}{2000}).
\newblock \bibinfo{title}{The unscented {K}alman filter for nonlinear
  estimation}.
\newblock In {\it \bibinfo{booktitle}{Proceedings of the IEEE 2000 Adaptive
  Systems for Signal Processing, Communications, and Control Symposium (Cat.
  No. 00EX373)}\/} (pp. \bibinfo{pages}{153--158}).
\newblock \bibinfo{organization}{IEEE}.
\bibitem[{Woodburn et~al.(2003)Woodburn, Carrico and
  Wright}]{woodburn2003estimation}
\bibinfo{author}{Woodburn, J.}, \bibinfo{author}{Carrico, J.}, and
  \bibinfo{author}{Wright, J.~R.} (\bibinfo{year}{2003}).
\newblock \bibinfo{title}{Estimation of instantaneous maneuvers using a fixed
  interval smoother}.
\newblock {\it \bibinfo{journal}{Advances in the Astronautical Sciences}\/},
  {\it \bibinfo{volume}{116}\/}, \bibinfo{pages}{243--260}.
\bibitem[{Xu et~al.(2019)Xu, Chen, Huang, Bai et~al.}]{xu2019collision}
\bibinfo{author}{Xu, Z.}, \bibinfo{author}{Chen, X.}, \bibinfo{author}{Huang,
  Y.}, \bibinfo{author}{Bai, Y.} et~al. (\bibinfo{year}{2019}).
\newblock \bibinfo{title}{Collision prediction and avoidance for satellite
  ultra-close relative motion with zonotope-based reachable sets}.
\newblock {\it \bibinfo{journal}{Proceedings of the Institution of Mechanical
  Engineers, Part G: Journal of Aerospace Engineering}\/},  {\it
  \bibinfo{volume}{233}\/}\bibinfo{issue}{(11)}, \bibinfo{pages}{3920--3937}.
\bibitem[{Ye et~al.(2021)Ye, Hua, Chuankai et~al.}]{ye2021maneuver}
\bibinfo{author}{Ye, L.}, \bibinfo{author}{Hua, Z.}, \bibinfo{author}{Chuankai,
  L.} et~al. (\bibinfo{year}{2021}).
\newblock \bibinfo{title}{Maneuver detection and tracking of a space target
  based on a joint filter model}.
\newblock {\it \bibinfo{journal}{Asian Journal of Control}\/},  {\it
  \bibinfo{volume}{23}\/}\bibinfo{issue}{(3)}, \bibinfo{pages}{1441--1453}.

\end{thebibliography}

\end{document}